\pdfoutput=1

\documentclass[a4paper]{article}
\usepackage{amsmath}
\usepackage{amsfonts,amssymb,amsthm}
\usepackage{color,enumerate}
\usepackage{hyperref}
\usepackage{graphicx}
\usepackage{spverbatim}
\usepackage{amsfonts}
\usepackage[total={6.3in, 10in}]{geometry}
\usepackage[normalem]{ulem}
\usepackage{siunitx,booktabs}
\usepackage[capitalise, noabbrev]{cleveref}
\usepackage{algpseudocode}
\usepackage{algorithm}
\usepackage{mathtools}

\usepackage{mathtools}
\usepackage{subcaption}
\definecolor{sixclassRdYlBu5}{rgb}{0.57,0.75,0.86}
\usepackage{float}
\usepackage{authblk}

\usepackage{multirow}

\usepackage{tikz}
\usetikzlibrary{calc}
\usepackage{pgfplots}
\usetikzlibrary{pgfplots.groupplots}
\usepackage{grffile}
\pgfplotsset{compat=newest}
\usetikzlibrary{plotmarks}

\newcommand{\setR}{\mathbb{R}}
\newcommand{\setN}{\mathbb{N}}

\newcommand{\bu}{\textbf{u}}
\newcommand{\bC}{\textbf{C}}
\newcommand{\bG}{\textbf{G}}
\newcommand{\bt}{\pmb{\theta}}
\newcommand{\bd}{\textbf{d}}

\newcommand{\bW}{\textbf{W}}

\newcommand{\J}{J}
\newcommand{\JML}{J_\text{ML}}

\newcommand{\transpose}{\mathsf{T}}

\newcommand{\ANNs}{DNN$_s$}
\newcommand{\ANNv}{DNN$_v$}

\newcommand{\EnOpt}{EnOpt}
\newcommand{\FOMEnOpt}{FOM-EnOpt}
\newcommand{\ROMEnOpt}{Adaptive-ML-EnOpt}
\newcommand{\ROMEnOpts}{AML-EnOpt$_s$}
\newcommand{\ROMEnOptv}{AML-EnOpt$_v$}

\renewcommand{\Return}[1]{\State{\textbf{return} #1}}

\theoremstyle{definition}
\newtheorem{definition}{Definition}

\definecolor{color0}{rgb}{0.65,0,0.15}
\definecolor{color1}{rgb}{0.84,0.19,0.15}
\definecolor{color2}{rgb}{0.96,0.43,0.26}
\definecolor{color3}{rgb}{0.99,0.68,0.38}
\definecolor{color4}{rgb}{1,0.88,0.56}
\definecolor{color5}{rgb}{0.67,0.85,0.91}
\definecolor{color6}{rgb}{0.27,0.46,0.71}

\def\FOMEnOptColor{color0}
\def\ROMEnOptsColor{color3}
\def\ROMEnOptvColor{color6}
\def\FOMValueColor{color1}

\def\FOMEnOptMarker{square*}
\def\ROMEnOptsMarker{triangle*}
\def\ROMEnOptvMarker{diamond*}
\def\FOMValueMarker{*}

\begin{document}

\title{Adaptive machine learning based surrogate modeling to accelerate PDE-constrained optimization in enhanced oil recovery}
\author[2]{Tim Keil}
\author[2]{Hendrik Kleikamp}
\author[1]{Rolf J. Lorentzen}
\author[1]{Micheal B. Oguntola\footnote{Corresponding author}}
\author[2]{\qquad Mario Ohlberger}

\affil[1]{NORCE-Norwegian Research Center AS, 5838, Bergen, Norway {\tt micheal.b.oguntola@uis.no, rolo@norceresearch.no}.}

\affil[2]{Institute for Analysis and Numerics, Mathematics Münster, University of Münster, Einsteinstrasse 62, 48149 Münster, Germany {\tt tim.keil@uni-muenster.de, hendrik.kleikamp@uni-muenster.de, mario.ohlberger@uni-muenster.de}.}

\maketitle

\begin{abstract}
\noindent In this contribution, we develop an efficient surrogate modeling framework for simulation-based optimization of enhanced oil recovery, where we particularly focus on 
polymer flooding. 
The computational approach is based on an adaptive training procedure of a neural network that directly approximates an input-output map of the underlying PDE-constrained optimization problem. The training process thereby focuses on the construction of an accurate surrogate model solely related to the optimization path of an outer iterative optimization loop. True evaluations of the objective function are used to finally obtain certified results. Numerical experiments are given to evaluate the accuracy and efficiency of the approach for a heterogeneous five-spot benchmark problem.
\end{abstract}

\section{Introduction}
Water flooding remains the most frequently used secondary oil recovery method. However, the percentage of original oil in place left after the cessation of water flooding in many reservoir fields is estimated to be as high as 50 - 70\% \cite{van2016well,pancholi2020experimental,zhang2017well}. The reduced performance of water flooding leading to the sizable leftover of oil has been linked to many factors such as the presence of unfavorable mobility ratios (due to heavy oil), high level of heterogeneity (in porosity and permeability), etc., in the reservoir \cite{gudina2015biosurfactant}. For these reasons, enhanced oil recovery (EOR) methods are employed to improve the performance of water flooding in order to increase oil production and minimize environmental stress.
\par
Polymer flooding is a matured chemical EOR method, suitable for heavy oil reservoir development, with over four decades of practical applications \cite{abidin2012polymers,wang2008key}. It involves injecting long chains of high-molecular-weight soluble polymers along with water flooding. The polymer EOR mechanism includes reducing mobility ratios of the oil-water system and early water breakthrough in the reservoir by increasing the viscosity of injected water and consequently improving vertical and aerial sweep efficiencies of the injected fluid. 
\par
The EOR process of polymer flooding can significantly increase the oil production \cite{wang2008key}. However, compared with water flooding, the operational cost and the risk associated with polymer flooding are higher. More so, since injecting more than necessary polymer into the reservoir can lead to insignificant oil increment, it is imperative to optimize the injection strategy of polymer flooding for field application to avoid unnecessarily high operational costs with no profit. 
\par
Conventionally, a reservoir simulation model is combined with a numerical optimization technique to determine an optimal control (including water rates, polymer concentrations of injection wells, liquid rates, or bottom hole pressures of production wells) for polymer flooding.
The aim is to maximize a given reservoir performance measure (RPM), such as the total oil production or the net present value (NPV) function over the reservoir life.
The simulation model is usually a complex numerical reservoir simulator that requires substantial data accounting for geology and geometry of the reservoir or rock and fluid properties.
In this study, the model simulates the oil reservoir response (inform of fluid production) to a given polymer flooding control per time.
On this account, we estimate the RPM of a given control strategy. 
\par
Further, the complexity of a reservoir simulator leads to a high computational effort for simulating a given polymer flooding scenario.
It contributes to the inefficiency of gradient-based solution techniques (e.g., the ensemble-based optimization (EnOpt) method) for polymer EOR optimization problems, since the (approximate) gradient of the objective functional with respect to the control variables requires several function evaluations, with each relying on a time-consuming polymer model simulation \cite{oguntola2020robust,xu2018production,zhou2013optimal}.
More so, for large-scale polymer problems discretized into a large number of grid cells, a single model evaluation may take several hours to complete.
For this reason, we propose a machine-learning-based approach to approximate the computationally demanding objective functional.
\par
In classical approaches of model order reduction or surrogate modeling, the expensive evaluation of the objective functional due to the PDE constraints is replaced by an a priori trained surrogate model that can be efficiently evaluated with respect to the optimization parameters.
In this work, however, we make use of an adaptive surrogate modeling approach, where a surrogate model is constructed during the outer optimization loop through adaptive learning that is targeted towards an accurate input-output map in the vicinity of the chosen parameters during the optimization loop.
The overall algorithm thus combines costly full order model (FOM) evaluations, training of machine learning (ML) based surrogate models, as well as evaluations of the successively trained ML models.
In model reduction for parameterized systems \cite{MR3701994}, such adaptive enrichment approaches have been recently proposed and successfully applied
in the context of PDE constrained parameter optimization, e.g., in combination with trust-region optimization \cite{Zahr2015,MR4269464,banholzer2020adaptive}.
Recently, in \cite{GHI+2021,HOS2021} first ideas were presented to combine online enrichment for reduced-order models (ROMs) with machine learning-based surrogate modeling.
In this contribution, we use feedforward deep neural networks (DNNs) to obtain surrogate models of the underlying input-output map that directly map the optimization parameters to the output of the objective functional.
\par
Artificial neural networks also gained attention in the context of
enhanced oil recovery in recent years, see \cite{hossein2021,cheraghi2021,ahmadi2015}, for instance.
However, these approaches mainly focus on accelerating the evaluation of the costly objective
function without providing a way to solve polymer EOR optimization problems
using the proposed surrogate models. In \cite{golzari2015}, the authors
describe an algorithm to obtain a global surrogate model that is applied as
a replacement for the objective functional in a genetic algorithm. The global
approximation of the objective is computed a priori before applying the
optimization routine. In \cite{lee2011}, artificial neural networks are
employed to facilitate the decision process for a specific EOR method.
\par
Concerning acceleration of PDE-constrained optimization in general, 
DNNs are, for instance, used in \cite{lye2021} to replace costly
simulations within the optimization loops by evaluations of surrogate
models. The main idea of the \emph{ISMO} algorithm described in
\cite{lye2021} is to run multiple parallel optimization routines
starting from different initial guesses and to construct DNN surrogate models
using training data collected at the final iterates of these optimization
algorithms. The training data is computed by costly evaluations of the
exact objective functional (involving the solution of PDEs). In contrast, the optimization
routines use the respective surrogate model to speed up the
computations. Iteratively, a surrogate model is built to
approximate the true objective functional near local optima. The
approximation quality also serves as the stopping criterion of the algorithm.
Another approach involving physics-informed deep operator networks 
to accelerate PDE-constrained optimization in a self-supervised manner 
has recently been suggested in~\cite{wang2021fast}. 
\par
The idea of not having a global surrogate model, but only approximations of
the objective functional that are locally accurate, is also one of the main
motivations for our algorithm.
In contrast to the procedure in \cite{lye2021} described previously, we
iteratively construct DNN surrogate models tailored towards the objective
function along a single optimization path. We consider only a single
initial guess but check for convergence by taking into account the true
objective functional. This stopping criterion certifies that
the resulting control is approximately a (local) optimum of the true objective functional
and not only of the surrogate. Further, we do not assume that the
derivative of the DNN surrogate with respect to its inputs is available
but reuse the EnOpt procedure when optimizing with the surrogate model.\\
\par

The remainder of this article is organized as follows. In \cref{sec:2}
we introduced the polymer flooding model for EOR and formulate 
an optimization problem for the economic value of the reservoir response. 
\cref{sec:3} introduces a classical ensemble based optimization algorithm 
based on a FOM approximation of the polymer flooding model.
Feedforward DNNs to approximate the input-output map are introduced in \cref{sec:4}. 
In \cref{sec:5}, we finally present and discuss our new adaptive FOM-ML-based optimization
algorithm, which is evaluated numerically for a five-spot benchmark problem in \cref{sec:6}. 
Last but not least, a conclusion and outlook is given in \cref{sec:7}.

\section{Optimization of polymer flooding in enhanced oil recovery}
\label{sec:2}
The problem of predicting the optimal injection strategy of the polymer EOR
method can be formulated as a constrained optimization problem. The setup
involves solving a maximization problem in which the objective function, the RPM, is defined on a given set of controllable 
variables. For the polymer EOR method, a complete set of control variables 
includes the concentration (and hence volume size of the polymer) and control 
variables (such as water injection rate, oil production rate, and/or bottom 
hole pressure for the injecting or producing wells) for water flooding over the 
producing lifespan of the reservoir. 

\subsection{Polymer flooding model}

As mentioned in the introduction, the optimization process is usually performed on a simulation model of the real reservoir~\cite{sarma2006efficient}. Here, we consider a polymer flooding simulation model, which is an extension of the black-oil model with a continuity equation for the polymer component \cite{zhou2013optimal,flo2019open}. The black-oil model is a special multi-component multi-phase flow model with no diffusion among the fluid components \cite{bao2017fully}. It assumes that all hydrocarbon species are considered as two components, namely, oil and gas at surface conditions, and can be partially or entirely dissolved in each other to form the oil and gas phases. Further, there is an aqueous phase that consists of only one component called water.
\par
For brevity, we first state the polymer flooding model without mentioning the dependence on the controls and geological parameters explicitly. These dependencies are described in more detail after depicting the model. Hence, in what follows, we assume that fixed sets of controls and geological parameters are given.
\par
In the polymer model, usually, it is assumed that polymer forms an additional component transported in the aqueous phase of the Black-oil model and has no effect on the oil phase. We identify those quantities associated with the water, oil, gas, and polymer components with subscripts W, O, G, and P. In general, the polymer model consists of the following system of partial differential equations:
\begin{subequations}\label{equ:polymer_model}
	\begin{align}
		\frac{\partial  }{\partial t}(\phi b_\text{W} s_\text{W}) + \nabla\cdot b_\text{W}\textbf{v}_\text{W} = q_\text{W},\\
		\frac{\partial }{\partial t}\phi (b_\text{O} s_\text{O}+ r_\text{OG}b_\text{G}s_\text{G}) + \nabla\cdot (b_\text{O}\textbf{v}_\text{O} + r_\text{OG}b_\text{G}\textbf{v}_\text{G}) = q_\text{O},\\
		\frac{\partial }{\partial t}\phi (b_\text{G} s_\text{G}+ r_\text{GO}b_\text{O}s_\text{O}) + \nabla\cdot (b_\text{G}\textbf{v}_\text{G} + r_\text{GO}b_\text{O}\textbf{v}_\text{O}) = q_\text{G},\\
		\frac{\partial }{\partial t}\Big[\phi(1 - s_{ipv})s_\text{W} + \frac{\rho_r c_a}{b_\text{W}c}(1 - \phi)\Big] + \nabla\cdot \textbf{v}_\text{P} = q_\text{W},
	\end{align}
\end{subequations}
where $\phi$ is the rock porosity, $s_{\alpha}, b_{\alpha}, q_{\alpha}$, and $\textbf{v}_{\alpha}$ denote the (unknown) saturation, inverse formation-volume factor (depending on the respective density $\rho_{\alpha}$), volumetric source (flow rate per unit volume), and Darcy's flux of phase $\alpha\in\{\text{W},\text{O},\text{G}\}$, and $r_\text{OG}$ and $r_\text{GO}$ denote the oil-gas and gas-oil ratios. The quantities $\textbf{v}_\text{P}, c_a, s_{ipv}$, and $c$ denote the Darcy's flux, adsorption concentration, inaccessible pore volume, and concentration of the polymer solution, and $\rho_r$ is the density of the reservoir rock.
\par
In addition to the system~\eqref{equ:polymer_model}, empirical closure equations for relative permeabilities and capillary pressure in three-phase flow in porous media are applied. Here, the unknown primary variables are phase saturations $s_{\alpha}$ (or component accumulations) and pressures $p_{\alpha}$, and thus, appropriate initial and boundary conditions are defined. 
\par
Based on the type of injection and/or production well (e.g., vertical, horizontal, or multi-segment), a suitable well model~\cite{holmes1998application,holmes1983enhancements} is coupled with~\eqref{equ:polymer_model} to measure the volumetric flow rates, which depend on the state of the reservoir. A standard well model for vertical wells is given as follows.
\par
The volumetric flow rates $q_{\alpha}$ for $\alpha\in\{\text{W}, \text{O}, \text{G}\}$ in a multi-phase polymer model are computed using a semi-analytical model according to~\cite{chen2007reservoir,holmes1998application} and are given by
\begin{subequations}\label{equ:well_equation}
	\begin{align}
		q_{\text{W}} &= \frac{k_\text{RW}(s_{\text{W}})}{\mu_{\text{W},\text{eff}}R_k(c)} WI(p_{\text{bh}} - p_{\text{W}} - \rho_{\text{W}}g(z_{\text{bh}}-z)),\\
		q_{\text{O}} &= \frac{k_\text{RO}(s_{\text{O}})}{\mu_{\text{O},\text{eff}}} WI(p_{\text{bh}} - p_{\text{O}} - \rho_{\text{O}}g(z_{\text{bh}}-z)),\\
		q_{\text{G}} &= \frac{k_\text{RG}(s_{\text{G}})}{\mu_{\text{G},\text{eff}}} WI(p_{\text{bh}} - p_{\text{G}} - \rho_{\text{G}}g(z_{\text{bh}}-z)).
	\end{align}
\end{subequations}
Here, $k_{\text{R}\alpha}(s_{\alpha})$, $\rho_{\alpha}$, $p_{\alpha}$, and $\mu_{\alpha,\text{eff}}$ are the saturation-dependent relative permeability,  density, pressure, and effective viscosity of phase $\alpha\in\{\text{W}, \text{O}, \text{G}\}$, $WI$ is the well index, $z_{\text{bh}}$ is the well datum level depth, $p_{\text{bh}}$ is the bottom hole pressure at the well datum level, $z$ is the depth, $R_k(c)$ models the reduced permeability experienced by the water-polymer mixture, and $g$ is the magnitude of the gravitational acceleration. 
\par
Individual wells are usually controlled by surface flow rates or bottom hole pressures. Additional equations which enforce limit values for the component rates and bottom-hole pressures are 
\begin{align*}
	p_{\text{bh}} - p_{\text{bh}}^{\text{limit}} &\leq 0,\\
	q_{\alpha} - q_{\alpha}^{\text{limit}} &\leq 0,
\end{align*}
where $q_{\alpha}^{\text{limit}}$ is the desired surface-volume rate limit for component $\alpha$, e.g., field oil rate at the production well, and $p_{\text{bh}}^{\text{limit}}$ is the desired bottom-hole pressure limit. Also, logic constraints to determine what happens if the computed rates or pressures violate the operational constraints, in which case a well may switch from rate control to pressure control, etc., are imposed.
\par
If $q_{\alpha,i}$ is the field volumetric flow rate (in sm$^3/$day) of component $\alpha\in\{\text{W},\text{O},\text{G}\}$ in the production wells over the time interval $\Delta t_i$, the field production total (in sm$^3$) of the component $\alpha$ is given as $Q_{\alpha\text{P},i} = q_{\alpha,i}\Delta t_i$. For polymer production total (in kg), $Q_{\text{PP},i} = c_{\text{L}}q_{\text{W},i}\Delta t_i$, where $c_\text{L}$ is the leftover field polymer concentration (in kg$/$sm$^3$) after adsorption. Injection quantities $Q_{\text{PI},i}$ and $Q_{\text{WI},i}$ are computed similarly, however with volumetric flow rates in the injection wells.
\par
As already mentioned above, the solution of the polymer flooding model
stated in~\eqref{equ:polymer_model} depends on a given control vector $\bu$,
see~\cref{subsec:opt_pb} for a detailed description of the components of the
control vector, and a set of geological properties $\bt$. Consequently, 
all involved unknowns depend on $\bu$ and $\bt$ and the same holds for $q_\text{W}$,
$q_\text{O}$, and $q_\text{G}$. From now on, we thus write
$Q_{\alpha\text{P},i}(\bu,\bt)$ for the field production total of component
$\alpha\in\{\text{W},\text{O},\text{G}\}$, depending on the controls $\bu$
and the parameters $\bt$, within the time interval $\Delta t_i$, similar as
above. We further write $Q_\text{PP}(\bu,\bt)$ for the polymer production
total, and $Q_{\text{PI},i}(\bu,\bt)$ and $Q_{\text{WI},i}(\bu,\bt)$ for the
polymer and water injection.

\subsection{Optimization of the economic value of the reservoir response}\label{subsec:opt_pb}

This study considers the annually discounted net present value (NPV) function as the RPM, similar to the one in~\cite{lu2020joint,oguntola2020robust}. The NPV function is related to the control variables through the polymer simulation model~\eqref{equ:polymer_model}. For every polymer control strategy, the NPV function evaluates the economic value of the reservoir response. Also, because the injection and production facilities have limited capacity, the control variables are subject to bound constraints.
\par
Suppose that the geological properties of the oil reservoir of interest, such as porosity, permeability, etc., are known
and denoted by $\bt$.
Let $\mathcal{D} = \setR^{N_u}$ be the domain of control vectors of polymer flooding for the given reservoir, such that
\[
	\bu = \left[u_1^1, u_2^1,\dots,u_{N_w}^{1},\dots,
	u_{1}^{N_t},u_{2}^{N_t},\dots,u_{N_w}^{N_t}\right]^\transpose,
\]
where $\transpose$ means transpose.
The subscript of each component of $\bu$ denotes the well index,
the superscript is the control time step, $N_w$ and $N_t$ denote the number of wells
and time steps for each well, respectively, and $N_u = N_w\cdot N_t$ is the total
number of control variables. Each component $u_j^i$ in $\bf u$ represents a
control type (e.g., polymer concentration or injection rate, oil or water rate,
bottom hole pressure) of well $j$ at the time step $i$. 
\par
The $N_u$-dimensional optimization problem for polymer flooding is to find the optimal
$\bu \in \mathcal{D}$ that maximizes the NPV function subject to bound constraints. That is
\begin{subequations} \label{equ:opt_problem}
	\begin{align}
		 &\mathop{\text{maximize}}_{\bu\in \mathcal{D}} ~\J(\bu, \bt) \coloneqq
		\sum_{i=1}^{N_t}\frac{\J_i(\bu,\bt)}{(1+d_{\tau})^{\frac{t_i}{\tau}}}\label{equ:fom-objective}\\
		 &\text{with}\nonumber\\
		 &\qquad \J_i(\bu,\bt) \coloneqq r_{\text{OP}}Q_{\text{OP},i}(\bu,\bt) +
		r_{\text{GP}}Q_{\text{GP},i}(\bu,\bt)-R_i(\bu,\bt),\nonumber\\
		 &\qquad R_i(\bu,\bt) \coloneqq r_{\text{WI}}Q_{\text{WI},i}(\bu,\bt) +
			r_{\text{WP}}Q_{\text{WP},i}(\bu,\bt) +r_{\text{PI}}Q_{\text{PI},i}(\bu,\bt) +
		r_{\text{PP}}Q_{\text{PP},i}(\bu,\bt),\nonumber\\
		&\text{subject to}\nonumber\\
		&\qquad u_j^{\text{low}} \leq u_j^i\leq u_j^{\text{upp}}
		\quad\text{for all } j=1,\dots,N_w,~i = 1,\dots,N_t, \label{equ:constraints}
	\end{align}
\end{subequations}
where $\J_i$ denotes the cumulative NPV value in the $i$-th simulation time
step. Further, $d_{\tau}$ is the discount rate for a period of $\tau$ days, $t_i$ is the cumulative time (in days) starting from the beginning of 
production up to the $i$-th time step, and $\Delta t_i \coloneqq t_i - t_{i-1}$ is the 
time difference (in days) between the time steps $t_i$ and $t_{i-1}$. The scalars
$r_{\text{OP}}, r_{\text{GP}}, r_{\text{WI}}$ and $r_{\text{WP}}$ denote the prices of oil and gas production 
and the cost of handling water injection and production (in USD/sm$^3$) 
respectively, and $r_{\text{PI}}$ and $r_{\text{PP}}$ are the costs of polymer injection and 
production (in USD/kg). In addition, $Q_{\text{WI},i}$ and $Q_{\text{PI},i}$ are the total 
water injection (in sm$^3$) and total polymer injection or slug 
size (in kg) over the time interval $\Delta t_i$. The quantities $Q_{\text{OP},i}$, 
$Q_{\text{WP},i}$ and $ Q_{\text{GP},i}$ denote the total oil, water and gas productions 
(in sm$^3$) over the time interval $\Delta t_i$, while $Q_{\text{PP},i}$ 
represents the total polymer production (in kg) over the time interval $\Delta t_i$.
The quantities $Q_{\text{OP},i}$, $Q_{\text{WI},i}$, $Q_{\text{WP},i}$, $Q_{\text{GP},i}$,  $Q_{\text{PI},i}$,
and $Q_{\text{PP},i}$ are computed at each control time 
step $i$ for given $\bu$ and fixed $\bt$ from the polymer flooding model~\eqref{equ:polymer_model} and the well equations~\eqref{equ:well_equation}.
\par
The evaluation of the objective function $\J$ in~\eqref{equ:fom-objective} shall be referred to as the full order model 
(FOM) function evaluation in the remainder of this study. Therefore, the 
constrained optimization problem presented in~\eqref{equ:opt_problem} can be interpreted as
the FOM optimization problem for polymer flooding, given a suitable discretization of the system~\eqref{equ:polymer_model} (see~\cref{sec:implementational-details} for details on the discretization).
Also, because~$\bt$ is fixed during the optimization process, $\J$ is considered a function of~$\bu$
only, and hence we often write~$\J(\bu)$ and $\J_i(\bu)$. The solution method utilized for this optimization
problem is presented in the next section.

\section{Ensemble based optimization algorithm}
\label{sec:FOM_enopt_algorithm}\label{sec:3}
In this work, the FOM solution to problem~\eqref{equ:opt_problem} follows
from the application of the adaptive ensemble-based optimization (\EnOpt)
method analogous to the one presented in~\cite{oguntola2020robust,
chen2009efficient,stordal2016theoretical}. We again emphasize that we restrict our attention to a fixed
choice of geological parameters $\bt$. Since we apply the \EnOpt\ algorithm
later on in our surrogate-based algorithm to a function different from
$\J$, we subsequently begin by describing the algorithm in its general
form. Afterwards, we discuss the application of the \EnOpt\ algorithm
to the objective function $\J$ and the resulting computational costs.

\subsection{Optimization algorithm for a general objective function}
\label{sec:EnOpt_general}
In what follows, we describe the \EnOpt\ algorithm for a general objective function $F\colon \setR^{N_u} \to \setR$ to iteratively solve the optimization problem
\begin{subequations} \label{equ:general_opt_problem}
	\begin{align}
		 &\mathop{\text{maximize}}_{\bu\in \mathcal{D}} ~F(\bu) \label{equ:general-fom-objective}\\
		&\text{subject to}
		\quad u_j^{\text{low}} \leq u_j^i\leq u_j^{\text{upp}}
		\quad\text{for all } j=1,\dots,N_w,~i = 1,\dots,N_t. \label{equ:general-constraints}
	\end{align}
\end{subequations}
The \EnOpt\ method is an iterative method in which one starts with an initial guess $\bu_0$ that is usually
based on experimental facts in such a way that the underlying constraints in~\eqref{equ:general-constraints} are satisfied.
We sequentially seek for an improved approximate solution $\bu$ that maximizes 
$F(\bu)$ using a preconditioned (with covariance matrix adaptation) gradient 
ascent method given by
\begin{align}\label{equ:UpdatingControl}
    \hat{\bu}_{k+1}&=\bu_k+\beta_{k}\,\bd_k,\\
    \bd_k&\approx\frac{\bC_{\bu_k}^k\bG_{k}}{\lVert\bC^k_{\bu_k}\bG_{k}\rVert_{\infty}},
\end{align}
where $k=0,1,2,\dots$ is the index of the optimization iteration. The tuning 
parameter $\beta_{k}$ for the step size is computed using an auxiliary line 
search~\cite{nocedal2006numerical} and is selected such that $0<\beta_k\leq 1$.
Furthermore, $\bC_{\bu_k}^k\in\setR^{N_u\times N_u}$ denotes the user-defined covariance 
matrix of the control variables at the $k$-th iteration and $\bG_k\in\setR^{N_u}$ is the 
approximate gradient of $F$ with respect to the control variables, preconditioned with
$\bC_{\bu_k}^k$ to obtain the search direction at the $k$-th iteration.
\par
To ensure that the constraints in~\eqref{equ:UpdatingControl} are satisfied, the original
solution domain of the control variables is projected to the set of admissible controls
$\mathcal{D}_\text{ad}$, defined as
\begin{align}
    \mathcal{D}_\text{ad} \coloneqq \{\bu\in\mathcal{D}: u_j^{\text{low}} \leq u_j^i\leq 
u_j^{\text{upp}}\text{ for all } j = 1,\dots,N_w, ~i = 1,\dots,N_t\},
\end{align}
which corresponds to the constraints in~\eqref{equ:general-constraints}.
The updating scheme in~\eqref{equ:UpdatingControl} is performed in $\mathcal{D}_\text{ad}$.
We utilize a component-wise projection $P_{\mathcal{D}_\text{ad}}\colon\mathcal{D}\to\mathcal{D}_\text{ad}$
on the update $\hat{\bu}_{k+1}\in\mathcal{D}$, such that
\begin{align}\label{equ:UpdatingControlWithProjectn}
    \bu_{k+1} &= P_{\mathcal{D}_\text{ad}}(\hat{\bu}_{k+1}) \in \mathcal{D}_\text{ad}.
\end{align}
In practical applications, it is not common to have controls at different wells 
to correlate, but the controls may vary smoothly with time at individual wells. 
Hence, the use of $\bC_{\bu_k}^k$ in Equation~\eqref{equ:UpdatingControl}
enforces this regularization on the control updates. At $k=0$, we utilize a 
temporal covariance function given by
\begin{align}\label{equ:CorrelationFunctn}
	\operatorname{Cov}\left(u_j^i,u_j^{i+h}\right) = 
\sigma^2_j\rho^h\left(\frac{1}{1-\rho^2}\right),\qquad\text{for all }h\in \{0,\dots,N_t-i\},
\end{align}
from a stationary auto regression of order 1 {(i.e., AR(1))}
model~\cite{montgomery2015introduction} to compute $\bC_{\bu_0}^0$ with an
assumption that controls of different wells are uncorrelated.
The variance for the well $j$ is given by $\sigma^2_j>0$, and $\rho\in(-1,1)$ is the correlation coefficient 
used to introduce a level of dependence  between controls of individual wells at 
different control time steps (since the AR(1) model is stationary).
\par
The formulation above gives rise to a block diagonal matrix $\bC_{\bu_0}^0$,
which is updated by matrices with rank one at subsequent iterations,
using the statistical method presented in~\cite{stordal2016theoretical}, to obtain an
improved covariance matrix $\bC_{\bu_k}^{k}$. For this reason, the solution 
method in Equation~\eqref{equ:UpdatingControl} is referred to as the adaptive \EnOpt\ algorithm.
\par
We compute the preconditioned approximate gradient 
$\bC_{\bu_k}^k\bG_{k}$ following the approach of the standard \EnOpt\ algorithm.
At the $k$-th iteration, we sample $N\in\setN$ control vectors 
$\bu_{k,m}\in\mathcal{D}_\text{ad}$,
for $m=1,\dots,N,$ from a multivariate Gaussian 
distribution with mean equal to the $k$-th control vector $\bu_k$ and 
covariance matrix given by $\bC_{\bu_k}^k$. Here, the additional subscript $m$ is 
used to differentiate the perturbed control vectors from the one obtained by 
Equation~\eqref{equ:UpdatingControl}.
The cross-covariance of the control vector $\bu_k$ and the objective function 
$F(\bu_k)$ at the $k$-th iteration is approximated according to~\cite{fonseca2017stochastic} as
\begin{align}\label{equ:Cross-covariance}
	\bC_{\bu_k,F}^k \coloneqq \frac{1}{N - 1}\sum_{m=1}^{N}(\bu_{k,m} - 
	\bu_k)\big(F(\bu_{k,m})- F(\bu_{k})\big).
\end{align}
Since $\bu_{k,m} \sim \mathcal{N}(\bu_k, \bC_{\bu_k}^k)$ for $m=1,\dots,N$, we 
assume in Equation~\eqref{equ:Cross-covariance} that the mean of
$\{\bu_{k,m}\}_{m=1}^{N}$ is approximated by $\bu_k$. By first-order Taylor
series expansion of $F$ about $\bu_k$, it can easily be deduced that
Equation~\eqref{equ:Cross-covariance} is an approximation of
$\bC_{\bu}^k\bG_{k}$ at the $k$-th iteration, that is
\begin{align}\label{equ:approximateCC}
	\bC_{\bu_k}^k\bG_{k} \approx \bC_{\bu_k,F}^k,
\end{align}
see~\cite{chen2009efficient, oguntola2021ensemble} for a detailed proof.
Therefore, we choose the search direction as $\bd_k=\bC_{\bu_k,F}^k/{\lVert\bC_{\bu_k,F}^k\rVert}_\infty$
in Equation~\eqref{equ:UpdatingControl}.
The updating scheme in Equation~\eqref{equ:UpdatingControl} is performed until the
convergence criterion
\begin{align}\label{equ:converCriteria}
	F(\bu_{k}) \leq F(\bu_{k-1}) + \varepsilon
\end{align}
is satisfied, where $\varepsilon>0$ is a specified tolerance.
\par
To conclude, for an arbitrary objective function $F$, the \EnOpt\ procedure described in this section is summarized in \cref{alg:EnOpt}.
In this algorithm, the \textproc{OptStep} function replicates a single optimization step in the \EnOpt\ procedure and is detailed in \cref{alg:OptStep}.
We note that returning the set of function values $T_{k+1}$ does not play a role in \cref{alg:EnOpt} but is crucial for training the surrogate model in \cref{sec:ann_enopt_algorithm}.
The line search procedure \textproc{LineSearch} can be found in \cref{alg:linesearch}.

\begin{algorithm}
	\caption{\EnOpt\ algorithm}\label{alg:EnOpt}
    \begin{algorithmic}[1]
        \Require{function $F\colon\setR^{N_u}\to\setR$ for which to solve~\eqref{equ:general_opt_problem};
                 initial guess $\bu_0\in\setR^{N_u}$,
                 sample size $N\in\setN$,
                 tolerance $\varepsilon > 0$,
                 maximum number of iterations $k^*$,
                 initial step size $\beta>0$,
                 step size contraction $r\in (0,1)$,
                 maximum number of step size trials $\nu^*\in\setN$}
		\Ensure{approximate solution $\bu^*\in\setR^{N_u}$ of~\eqref{equ:general_opt_problem}}
        \Function{\EnOpt[$F$]}{$\bu_0$, $N$, $\varepsilon$, $k^*$, $\beta$, $r$, $\nu^*$}
        \State{$\bu_1,\ T_1\leftarrow$ \Call{OptStep[$F$]}{$\bu_0$, $N$, $0$, $\beta$, $r$, $\nu^*$}\label{iteration_0}}
        \State{$k\leftarrow 1$}
        \While{$F(\bu_{k}) > F(\bu_{k-1}) + \varepsilon$ and $k < k^*$ \label{alg:termination}}
            \State{$\bu_{k+1},\ T_{k+1}\leftarrow$ \Call{OptStep[$F$]}{$\bu_k$, $N$, $k$, $\beta$, $r$, $\nu^*$}}
            \State{$k\leftarrow k+1$}
        \EndWhile
        \Return{$\bu^*\leftarrow\bu_{k}$}
        \EndFunction
	\end{algorithmic}
\end{algorithm}

\begin{algorithm}
	\caption{OptStep algorithm}\label{alg:OptStep}
	\begin{algorithmic}[1]
        \Require{function $F\colon\setR^{N_u}\to\setR$;
        		 current control vector $\bu_k\in\setR^{N_u}$,
				 sample size $N\in\setN$,
				 number of iteration $k$,
				 initial step size $\beta>0$,
				 step size contraction $r\in (0,1)$,
				 maximum number of step size trials $\nu^*\in\setN$}
        \Ensure{update $\bu_{k+1}\in\setR^{N_u}$ of the controls, set $T_{k+1}$ of $N$ pairs of the form $(\bu,F(\bu))$}
        \Function{OptStep[$F$]}{$\bu_k$, $N$, $k$, $\beta$, $r$, $\nu^*$}
        \If{$k = 0$}
			\State{Compute the initial covariance matrix $\bC^0_{\bu_0}$
				   using~\eqref{equ:CorrelationFunctn}}
        \Else
			\State{Compute the covariance matrix $\bC^k_{\bu_k}$
				   using the formulation in~\cite{stordal2016theoretical}}
        \EndIf
        \State{Sample $N$ control vectors $\{\bu_{k,j}\}_{j=1}^N$ from a
			   distribution $\mathcal{N}(\bu_{k}, \bC^k_{\bu_k})$}
        \State{Compute vector $\bC_{\bu_k,F}^k$ according to~\eqref{equ:Cross-covariance} and store values $\{F(\bu_{k,j})\}_{j=1}^N$\label{comp_samples}}
        \State{Compute the search direction
        	   $\bd_k=\bC_{\bu_k,F}^k/\lVert\bC_{\bu_k,F}^k\rVert_\infty$ \label{comp_dir}}
        \State{$\bu_{k+1}\leftarrow$ \Call{LineSearch[$F$]}{$\bu_k$, $\bd_k$, $\beta$, $r$, $\nu^*$}\label{comp_line}}
		\State{$T_{k+1}\leftarrow\{(\bu_{k,j},F(\bu_{k,j}))\}_{j=1}^N$ \label{training_data}}
        \Return{$\bu_{k+1}$, $T_{k+1}$}
        \EndFunction
	\end{algorithmic}
\end{algorithm}

\begin{algorithm}[H]
	\caption{Line search}\label{alg:linesearch}
	\begin{algorithmic}[1]
		\Require{function $F\colon\setR^{N_u}\to\setR$;
			current controls $\bu_k\in\setR^{N_u}$,
			search direction $\bd_k\in\setR^{N_u}$,
			initial step size $\beta>0$,
			step size contraction $r\in (0,1)$,
			maximum number of step size trials $\nu^*\in\setN$,
			tolerance $\varepsilon > 0$}
		\Ensure{update $\bu_{k+1}\in\setR^{N_u}$ of the controls}
		\Function{LineSearch[$F$]}{$\bu_k$, $\bd_k$, $\beta$, $r$, $\nu^*$}
		\State{$\beta_k\leftarrow\beta$}
		\State{Compute $\bu_{k+1}$ according to~\eqref{equ:UpdatingControlWithProjectn}}
		\State{$\nu\leftarrow 0$}
		\While{$F(\bu_{k+1}) - F(\bu_k) \leq \varepsilon$ and $\nu<\nu^*$}
		\State{$\beta_k\leftarrow r\,\beta_k$}
		\State{Compute $\bu_{k+1}$ according to~\eqref{equ:UpdatingControlWithProjectn}}
		\State{$\nu\leftarrow\nu + 1$}
		\EndWhile
		\State{\textbf{return} $\bu_{k+1}$}
		\EndFunction
	\end{algorithmic}
\end{algorithm}

\subsection{\FOMEnOpt\ algorithm for enhanced oil recovery}
Eventually, we are interested in solving the optimization problem~\eqref{equ:opt_problem} for polymer flooding in enhanced oil recovery. As already discussed in the introduction, our contribution is concerned with the development of a surrogate-based algorithm to reduce the computational costs for solving~\eqref{equ:opt_problem}. To this end, if the \EnOpt\ algorithm is used to maximize the function $\J$, defined in \cref{equ:fom-objective}, we refer to \cref{alg:EnOpt} as the \emph{\FOMEnOpt\ algorithm}. That is, the \FOMEnOpt\ algorithm is given as \textproc{\EnOpt[$\J$]}, see \cref{alg:FOMEnOpt}.

 \begin{algorithm}
 	\caption{\FOMEnOpt\ algorithm}\label{alg:FOMEnOpt}
 	\begin{algorithmic}[1]
 		\Require{initial guess $\bu_0\in\setR^{N_u}$,
 				 sample size $N\in\setN$,
 				 tolerance $\varepsilon > 0$,
 				 maximum number of iterations $k^*$,
 				 initial step size $\beta>0$,
 				 step size contraction $r\in (0,1)$,
 				 maximum number of step size trials $\nu^*\in\setN$}
 		\Ensure{approximate solution $\bu^*\in\setR^{N_u}$ of~\eqref{equ:opt_problem}}
 		\Function{\FOMEnOpt}{$\bu_0$, $N$, $\varepsilon$, $k^*$, $\beta$, $r$, $\nu^*$}
 			\Return{\Call{EnOpt[$\J$]}{$\bu_0$, $N$, $\varepsilon$, $k^*$ $\beta$, $r$, $\nu^*$}}
 		\EndFunction
 	\end{algorithmic}
 \end{algorithm}

As already indicated, we are concerned with the computational effort of the \FOMEnOpt\ algorithm.
Let us recall that evaluating $\J$ as in~\eqref{equ:fom-objective} has the complexity of the high-fidelity reservoir simulator, which, in itself, requires the solution of the discretized polymer flooding model equations~\eqref{equ:polymer_model}.
In \cref{alg:EnOpt}, the most expensive part is to call \textproc{OptStep[$\J$]}, which requires $N$ evaluations of $\J$ in Line~\ref{comp_samples} of \cref{alg:OptStep} such that the direction $\bd_k$ can be computed in Line~\ref{comp_dir}.
Furthermore, the line search in Line~\ref{comp_line} evaluates~$\J$ for every search step.
Suppose the simulation time for computing $\J$ is particularly large. In that case, the \FOMEnOpt\ algorithm can be extremely costly, especially if many optimization steps are required since \textproc{OptStep[$\J$]} is called at every iteration step.
In this case, all steps in \cref{alg:EnOpt} and \cref{alg:OptStep} that do not require evaluating $\J$ are computationally negligible.
\par
Since expensive FOM evaluations are very likely to happen for the presented
application, we aim to derive a surrogate-based algorithm that uses an approximation of $\J$ whenever possible and thus tries to reduce the number of calls of \textproc{OptStep[$\J$]}. Instead, FOM information is
reused whenever possible and only computed when necessary.
The following section introduces a machine-learning-based way for deriving suitable non-intrusive surrogate models.

\section{Neural networks as surrogate model for the input-output map}
\label{sec:anns}\label{sec:4}
Deep neural networks (DNNs) are machine learning algorithms suitable for approximating
functions without knowing their exact structure. Instead, DNNs can be fitted to approximately
reproduce known target values for a set of given inputs. Since DNNs learn from examples of labeled data,
they can be seen as \emph{supervised} learning algorithms. In contrast, \emph{unsupervised}
machine-learning algorithms try to detect hidden structures within unlabeled data.
See~\cite{hastie2009} for an exhaustive overview of supervised and unsupervised learning algorithms.
\par
A particular class of DNNs are \emph{feedforward neural networks}, in which no cyclic flow of information is allowed.
This study considers feedforward neural networks consisting of
(fully-connected) linear layers combined with a nonlinear activation function.
Our description of these types of DNNs is based on formal definitions that
can be found in~\cite{petersen2018} and~\cite{elbraechter2018}, for instance.
\par
Feedforward neural networks are used to approximate a given function $f\colon\setR^{N_\text{in}}\to\setR^{N_\text{out}}$
for a certain input dimension $N_\text{in}\in\setN$ and an output
dimension $N_\text{out}\in\setN$. To this end, let $L\in\setN$ denote the
\emph{number of layers} in the neural network,
and $N_\text{in}=N_0,N_1,\dots,N_{L-1},N_L=N_\text{out}\in\setN$ the
\emph{numbers of neurons} in each layer.
Furthermore, the \emph{weights} and \emph{biases} in layer $i\in\{1,\dots,L\}$ are denoted by
$W_i\in\setR^{N_i\times N_{i-1}}$ and $b_i\in\setR^{N_i}$. We assemble the weights and biases in an
$L$-tuple $\bW=\big((W_1,b_1),\dots,(W_L,b_L)\big)$. Moreover, let $\rho\colon\setR\to\setR$ be the
so-called \emph{activation function} and $\rho_n^*\colon \setR^n \to \setR^n$ the component-wise application of the activation
function $\rho$ for dimension $n\in\setN$, that is $\rho_n^*(y)\coloneqq
\left[\rho(y_1),\dots,\rho(y_n)\right]^\transpose\in\setR^n$ for~$y\in\setR^n$.
Then we can define the corresponding feedforward neural network in the following way:
\begin{definition}[Feedforward neural network]\label{def:DNN}
    The \emph{feedforward neural network} with weights and biases~$\bW$ and activation function~$\rho$
    for approximating $f\colon\setR^{N_\text{in}}\to\setR^{N_\text{out}}$, is defined as the function
    $\Phi_\bW\colon \setR^{N_\text{in}}\to\setR^{N_\text{out}}$.
    For a given input $x\in\setR^{N_\text{in}}$, the result $\Phi_\bW(x)\in\setR^{N_\text{out}}$ is computed as
    \begin{align*}
        \Phi_\bW(x) \coloneqq r_L(x),
    \end{align*}
    where $r_L\colon\setR^{N_\text{in}}\to\setR^{N_\text{out}}$ is defined in a recursive
    manner using the functions $r_i\colon\setR^{N_\text{in}}\to\setR^{N_i}$ for
    $i=0,\dots,L-1$, which are given by
    \begin{align*}
        r_L(x) &\coloneqq W_L\,r_{L-1}(x) + b_L,\\
        r_i(x) &\coloneqq \rho_{N_i}^*\left(W_{i}\,r_{i-1}(x) + b_i\right)&\text{for }i=1,\dots,L-1,\\
        r_0(x) &\coloneqq x.
    \end{align*}
\end{definition}

Fitting neural network weights and biases to a given function $f$ is accomplished by creating a
sample set $T_\text{train}=\{(x_1,f(x_1)),\dots,(x_n,f(x_n))\}
\subset X\times\setR^{N_\text{out}}$ (the so-called \emph{training set}),
consisting of inputs $x_i\in X$ from an input set $X\subset\setR^{N_\text{in}}$
and corresponding outputs $f(x_i)\in\setR^{N_\text{out}}$. The process of finding the weights $\bW$
such that $\Phi_\bW(x_i)\approx f(x_i)$ for $i=1,\dots,n$ is called \emph{training} of the
neural network. During the training, the weights and biases of the neural network $\Phi_\bW$
are iteratively adjusted such that a \emph{loss function}, which measures the
deviation of the output $\Phi_\bW(x_i)$ for a given input $x_i$ from the desired result $f(x_i)$, is minimized.
A common choice for the loss function is the \emph{mean squared error loss}
$\mathcal{L}\left(\Phi_\bW,T_\text{train}\right)$ given as
\begin{align}
    \mathcal{L}\left(\Phi_\bW,T_\text{train}\right) \coloneqq \sum\limits_{(x,y)\in T_\text{train}} \lVert\Phi_\bW(x) - y\rVert_2^2.
\end{align}
For a fixed \emph{architecture}, i.e.\ fixed number of layers $L$ and numbers of neurons
$N_0,\dots,N_L$ in each layer, we define the set of possible weights and biases $\Psi$ as
\begin{align*}
    \Psi\coloneqq \bigtimes_{i=1}^L \left(\setR^{N_i\times N_{i-1}}\times\setR^{N_i}\right).
\end{align*}
The set $\Psi$ contains $L$-tuples such that the matrices and vectors in each tuple
have suitable dimensions.
The aim of neural network training is to find weights and biases $\bW^*\in\Psi$ such that the
corresponding function $\Phi_{\bW^*}$ minimizes the loss function $\mathcal{L}$, i.e.
\begin{align}\label{equ:DNNOptimizationProblem}
    \bW^* = \arg\min\limits_{\bW\in\Psi}\ \mathcal{L}\left(\Phi_\bW,T_\text{train}\right).
\end{align}
There are several suitable optimization algorithms to approximate the solution
of~\eqref{equ:DNNOptimizationProblem} numerically. All of these methods require access to the gradient
of the loss function $\mathcal{L}$ with respect to the weights $\bW$ of the DNN, which can be computed
efficiently using an algorithm called \emph{backpropagation}, see~\cite{rumelhart1986}.
Popular examples of optimization algorithms used in neural network training are variants of
\emph{(stochastic) gradient descent methods}, see~\cite{bottou2018} for an overview. For small neural
networks with only a few layers and neurons, it is also possible to apply methods that use
or approximate higher-order derivatives of the loss function, for instance, the
\emph{L-BFGS} optimizer~\cite{liu1989}, which is a limited-memory variant of the \emph{BFGS}
method, see for instance Section~6.1 in~\cite{nocedal2006numerical}.
In the context of neural network training, each iteration of the optimizer is called \emph{epoch}.
Typically, a maximal number of epochs is prescribed for the optimizer to perform.
\par
To prevent a neural network from \emph{overfitting} the training data, we employ \emph{early stopping}~\cite{prechelt1997}.
In this method, the loss function is evaluated on a \emph{validation set} $T_\text{val}\subset X\times\setR^{N_\text{out}}$
after each epoch. The validation set is usually chosen to be disjoint from the training set, i.e.\
$T_\text{val}\cap T_\text{train}=\emptyset$. Let $\bW_k\in\Psi$ denote the weights in epoch $k\in\setN$.
In each epoch, the value $\mathcal{L}(\Phi_{\bW_k},T_\text{val})$ is computed, and if this value does
not decrease anymore over a prescribed number of consecutive epochs, the training is aborted.
This method ensures that the resulting neural network can perform well on unseen data
(that is assumed to have the same structure as the training data).
\par
The result of the optimization routine typically depends strongly on the initial values $\bW_0\in\Psi$
of the weights. There are several methods for initializing the weights of neural networks,
for instance, the so-called Kaiming initialization, see~\cite{he2015} for more details.
We perform \emph{multiple restarts} of the training algorithm using different initial values for the weights
to minimize the dependence of the resulting neural network on the weight initialization. Finally,
we select the neural network $\Phi_{\bW^*}$ that produced the smallest loss $\mathcal{L}(\Phi_{\bW^*},T_\text{train})
+ \mathcal{L}(\Phi_{\bW^*},T_\text{val})$ over all training restarts, i.e.\ the smallest combined
loss on the training and the validation set.
\par
Finding an appropriate neural network architecture can be difficult in practical applications.
Especially the number of layers and the number of neurons significantly influence the
approximation capabilities of the resulting neural network.
We call a layer \emph{hidden} if it is not an input or an output layer.
Neural networks with more than one hidden layer are called \emph{deep neural networks}.
See~\cite{yarotsky2017} for proofs that DNNs have an increased expressiveness.
In addition, there are lots of different activation functions available.
Typical examples include the \emph{rectified linear unit (ReLU)} $\rho(x)=\max(x,0)$, which is nowadays the
most popular activation function~\cite{lecun2015}, or the hyperbolic tangent $\rho(x)=\tanh(x)=\frac{e^{2x}-1}{e^{2x}+1}$.

\section{\ROMEnOpt\ algorithm using deep neural networks}
\label{sec:ann_enopt_algorithm}\label{sec:5}
The primary purpose of this work is to propose an adaptive machine-learning-based algorithm for avoiding expensive FOM evaluations as often as possible.
To this end, we first discuss the usage of DNNs for the NPV value and subsequently introduce the \ROMEnOpt\ algorithm.
   
\subsection{Surrogate models for the net present value}
As discussed in \cref{sec:FOM_enopt_algorithm}, we use DNNs to construct a surrogate model for the FOM objective functional~$\J$.
DNNs are particularly well suited for non-intrusive model reduction if the simulator is considered a black box with no direct access to solutions of the underlying PDEs.
In fact, given the formulation of the objective functional~\eqref{equ:fom-objective}, we assume to only have access to the respective components~$\J_i(\bu)$.
\par
Following the definition of a DNN in \cref{sec:anns}, two input-output maps can be used to approximate~$\J$.
We refer to the \emph{scalar-valued output} by considering $\J\colon \setR^{N_u} \to \setR$ as the input-output map.
Furthermore, we refer to the \emph{vector-valued output} if we make different use of the structure of $\J$ by writing $\J(\bu)=\delta^\transpose j(\bu)$ with
\begin{align*}
	j&\colon \setR^{N_u} \to \setR^{N_t}, \\
	j(\bu) &\coloneqq \left[\J_i(\bu)\right]_{i=1}^{N_t},
\end{align*}
and the vector $\delta\in\setR^{N_t}$, which includes the discount factors, is defined as
\begin{align*}
	\delta\coloneqq\left[\frac{1}{(1+d_{\tau})^{\frac{t_i}{\tau}}}\right]_{i=1}^{N_t}.
\end{align*}
In the scalar-valued case (\ANNs-approach), we directly construct a DNN for
$\J$ with a corresponding function $\Phi_{\bW_s}\colon \setR^{N_u} \to \setR$,
i.e.~we use a DNN with $N_\text{in}=N_u$ and $N_\text{out}=1$.
Instead, in the vector-valued case (\ANNv-approach), we construct a DNN for approximating $j$ with
a corresponding function $\Phi_{\bW_v}\colon \setR^{N_u} \to \setR^{N_t}$ and,
by using $\delta$, we indirectly approximate $\J$. This means that we
apply a DNN with input- and output-dimensions given by $N_\text{in}=N_u$ and
$N_\text{out}=N_t$, and multiply the result by $\delta$ whenever the respective
DNN is used for approximating $\J$.
The algorithm described below works for both cases, the scalar-valued and the vector-valued output.
Therefore, if access to the individual components of the vector-valued function $j$ is available, it is possible to run the algorithm with both versions.
The different neural network output sizes, and therefore, the various structures of the training data, might improve the DNN training results.
In our numerical experiment, we observe that the vector-valued DNN yields slightly better results than the scalar-valued DNN (see \cref{sec:numerical_experiments}).
Nevertheless, we consider both the scalar- and vector-valued approaches to discuss the case where the black box reservoir simulator produces only scalar-valued outputs.
\par
By the \ANNs- and \ANNv-approach, we thus construct a surrogate for the objective function for the optimization problem~\eqref{equ:opt_problem}.
It remains to explain a suitable and robust \EnOpt\ algorithm that takes advantage of a DNN but shows a similar convergence behavior as the FOM algorithm.
A common strategy is to construct a sufficiently accurate surrogate $\JML \in \{ \Phi_{\bW_s}, \delta^\transpose \Phi_{\bW_v} \}$ for the entire input space in a large offline time.
Following the \FOMEnOpt\ procedure from \cref{sec:FOM_enopt_algorithm}, given $\JML$, a surrogate-based procedure would then mean to set $F \coloneqq \JML$ in \cref{alg:EnOpt}.
However, no FOM stopping criterion would be used, and since no error control for the surrogate model is given, no certification of the surrogate-based procedure would be available.
Importantly, we remark that the input dimension $N_u$ of both DNN approaches is proportional to the number of time steps $N_t$ and the number of physical variables in the model~$N_w$.
Thus, dependent on the complexity of the reservoir simulation, $N_u$ may be large.
Consequently, it may not be possible to construct a surrogate model with a DNN that is accurate for the entire input space.
Even if it were possible to construct such a DNN, we would require prohibitively costly training for computing the training set, validation set, and weights.

\subsection{Adaptive algorithm}
To circumvent the issue of constructing a globally accurate surrogate, in what follows, we describe the adaptive machine learning \EnOpt\ algorithm (\ROMEnOpt).
In this algorithm, we incorporate the construction of the DNN into an outer optimization loop trained and certified by FOM quantities.
With respect to the \FOMEnOpt\ procedure, we remark that each FOM optimization step requires
$N$ evaluations of $\J$ for computing $\bd_k$.
To obtain an appropriately accurate direction, it is required that $N$ is chosen sufficiently large \cite{fonseca2015quantification}.
For the \ROMEnOpt\ procedure, we only use a single FOM-based optimization step at each outer iteration $k$.
Then, we use the $N$ evaluations of the FOM as data points for training a locally accurate surrogate~$\JML^k$.
Instead of proceeding with the FOM functional $\J$, we utilize the DNN to start an inner \EnOpt\ algorithm with $F=\JML^k$ as objective function in \cref{sec:EnOpt_general} and $\bu_k$ as initial guess.
Denote by $\bu_k^{(l)}$ the iterates of the inner optimization loop in the $k$-th outer iteration, i.e., in particular, we have~$\bu_k^{(0)}=\bu_k$.
According to~\eqref{equ:converCriteria}, the inner \EnOpt\ iteration terminates if the surrogate-based criterion
\begin{align}\label{equ:converCriteria_ML}
	\JML^k(\bu_k^{(l)}) \leq \JML^k(\bu_k^{(l-1)}) + \varepsilon_i
\end{align}
is met for a suitable tolerance $\varepsilon_i > 0$.
If the inner iteration terminates after $L$ iterations with a control~$\bu_k^{(L)}$, the next outer iterate $\bu_{k+1}$ is defined as $\bu_{k+1} := \bu_k^{(L)}$.
For a certified FOM-based stopping criterion of the outer optimization loop, given the iterate~$\bu_k$, we check whether the \FOMEnOpt\ procedure would, indeed, also stop at the same control point.
Thus, we perform a single FOM-based optimization step, which includes the computation of $\bd_k$ and the line search, and results in a control $\tilde{\bu}_k$.
For verifying whether the FOM optimization step successfully finds a sufficiently increasing point at outer iteration $k$, we consider the FOM termination criterion
\begin{align}\label{equ:converCriteria_Adaptive_ML}
J(\tilde{\bu}_{k}) \leq J(\bu_{k}) + \varepsilon_o,
\end{align}
where $\varepsilon_o > 0$ is a suitable tolerance.
If~\eqref{equ:converCriteria_Adaptive_ML} is fulfilled, no improvement of the objective function value using FOM optimization steps can be expected, and therefore we also terminate the \ROMEnOpt\ algorithm.
If instead,~\eqref{equ:converCriteria_Adaptive_ML} is not met, we use the computed training data (collected while computing~$\bd_k$) to retrain the DNN and restart an inner DNN-based \EnOpt\ algorithm.
We emphasize that the fully FOM-based stopping criterion constitutes a significant difference to what is proposed in \cite{lye2021}, where the termination criterion is based on the approximation quality of the surrogate model at the current iterate.
However, we saw in our experiments that such an approximation-based criterion might lead to an undesired early stopping of the algorithm.
\par
One may be concerned about the fact that the surrogate-based inner optimization routine produces a decreasing or stationary point.
For this reason, after every outer iteration $k$ of the \ROMEnOpt\ procedure, the inner DNN-optimization is only accepted after a sufficient increase, i.e.
\begin{equation}\label{eq:acceptance}
	\J(\bu_{k+1}) > \J(\bu_{k}) + \varepsilon_o.
\end{equation}
If an iterate is not accepted, we abort the algorithm.
Instead of aborting, one may proceed with an intermediate FOM optimization step.
We would further like to emphasize that the fulfillment of~\eqref{eq:acceptance}
also depends on the successful construction of the neural network, meaning that the parameters for the neural network are chosen appropriately.
If, instead,~\eqref{eq:acceptance} fails due to an inaccurate neural network,
an automatic variation of the parameters could be enforced to the neural network training, and the corresponding outer iteration should be repeated.
However, for the sake of simplicity and because it did not show any relevance in our numerical experiments, we do not specify approaches for the case that $\bu_{k+1}$ is not accepted due to~\eqref{eq:acceptance}.
\par
Regarding the choice of the different tolerances~$\varepsilon_i$
and~$\varepsilon_o$ for the inner and outer stopping criteria in the
\ROMEnOpt\ algorithm, we propose to choose a small value for~$\varepsilon_i$
similar to the tolerance~$\varepsilon$ in the \FOMEnOpt\ procedure.
The inner iterations are much cheaper due to the application of a fast
surrogate, such that a more significant amount of inner iterations is acceptable.
In contrast, we recommend selecting a larger tolerance~$\varepsilon_o$ to
perform fewer outer iterations for obtaining a considerable speed-up.
However, if maximum convergence w.r.t. the \FOMEnOpt\ algorithm is desired,
$\varepsilon_o$ is to be set equal to~$\varepsilon$.
\par
The above-explained \ROMEnOpt\ procedure is summarized in \cref{alg:ROM_EnOpt}.

\begin{algorithm}[H]
	\caption{\ROMEnOpt\ algorithm}\label{alg:ROM_EnOpt}
    \begin{algorithmic}[1]
        \Require{initial guess $\bu_0\in\setR^{N_u}$,
				 sample size $N\in\setN$,
                 tolerance $\varepsilon_o > 0$ for outer iterations,
                 tolerance $\varepsilon_i > 0$ for inner iterations,
                 maximum number of outer iterations $k_o^*$,
                 maximum number of inner iterations $k_i^*$,
				 DNN construction strategy $\text{CS}\in\{\text{\ANNs}, \text{\ANNv}\}$,
				 set of DNN-specific variables $V_\text{DNN}$ as discussed in~\cref{sec:anns} (e.g.~network architecture, loss function, training parameters),
                 initial step size $\beta>0$,
                 step size contraction $r\in (0,1)$,
                 maximum number of step size trials $\nu^*\in\setN$}
		\Ensure{approximate solution $\bu^*\in\setR^{N_u}$ of~\eqref{equ:opt_problem}}
		\Function{RomEnOpt}{$\bu_0$, $N$, $\varepsilon_o$, $\varepsilon_i$, $k_o^*$, $k_i^*$, $\text{CS}$, $V_\text{DNN}$, $\beta$, $r$, $\nu^*$}
		\State{$\tilde{\bu}_0,\ T_0\leftarrow$ \Call{OptStep[$\J$]}{$\bu_0$, $N$, $0$, $\beta$, $r$, $\nu^*$} \label{opt_step_1}}
		\State{$k\leftarrow 0$}
        \While{$\J(\tilde{\bu}_{k}) > \J(\bu_{k}) + \varepsilon_o$ and $k < k_o^*$ \label{alg:stopping_ROM}}
			\State{$\JML^k\leftarrow$ \Call{Train}{$T_k$, $\text{CS}$, $V_\text{DNN}$} \label{alg:training}}
			\State{$\bu_{k+1}\leftarrow$  \Call{EnOpt[$\JML^{k}$]}{$\bu_{k}$, $N$, $\varepsilon_i$, $k_i^*$ $\beta$, $r$, $\nu^*$}}
			\If{$\J(\bu_{k+1}) \leq \J(\bu_{k}) + \varepsilon_o$ \label{acceptance}}
				\Return{$\bu^*\leftarrow\bu_{k}$}
			\EndIf
			\State{$\tilde{\bu}_{k+1},\ T_{k+1}\leftarrow$ \Call{OptStep[$\J$]}{$\bu_{k+1}$, $N$, $k$, $\beta$, $r$, $\nu^*$} \label{opt_step_2}}
			\State{$k\leftarrow k+1$}
        \EndWhile
        \State{\textbf{return} $\bu^*\leftarrow \bu_{k}$}
        \EndFunction
	\end{algorithmic}
\end{algorithm}

The \textproc{Train} function performs the neural network training
procedure as described in \cref{sec:anns} and returns, depending on the
chosen DNN construction strategy, a function $\Phi_{\bW_s}$ or
$\Phi_{\bW_v}$ that approximates the FOM objective function $\J$. Particularly, the result of \textproc{Train} can be used as the function $F$ in the \textproc{EnOpt} procedure.
\par
The outer acceptance criterion~\eqref{eq:acceptance} is checked in Line~\ref{acceptance}. Using the FOM-based stopping criterion in Line~\ref{alg:stopping_ROM}, we ensure that the \ROMEnOpt\ algorithm has an equivalent stopping procedure as the \FOMEnOpt\ algorithm, see Line~\ref{alg:termination} in \cref{alg:EnOpt}.
However, the algorithm might terminate at a different (local) optimal point, which we also observe in the numerical experiments.
\par
Compared to the \FOMEnOpt\ procedure, we emphasize that, in the \ROMEnOpt\ algorithm, mainly the single calls of \textproc{OptStep[$\J$]} in Lines~\ref{opt_step_1} and~\ref{opt_step_2} have FOM complexity, scaling with the number of samples~$N$. Furthermore, the outer stopping criterion in Line~\ref{alg:stopping_ROM} and the conditions for acceptance in Line~\ref{acceptance} require a single FOM evaluation.
The construction of the surrogate makes use of FOM data that is already available from the FOM optimization steps in Lines~\ref{opt_step_1} and~\ref{opt_step_2}.
In addition, while the training data for Line~\ref{alg:training} is available from calling \textproc{OptStep[$\J$]}, the training function \textproc{Train} itself is relatively cheap.
Furthermore, calling \textproc{EnOpt[$\JML^{k}$]} has low computational effort since evaluating the surrogate~$\JML^{k}$ for a given control (i.e., performing a single forward pass through the neural network) is much faster than evaluating $\J$.
The primary motivation for the \ROMEnOpt\ algorithm is the idea that many of the costly FOM optimization steps in the \FOMEnOpt\ algorithm can be replaced by sequences of cheap calls of \textproc{EnOpt[$\JML^{k}$]} with the surrogate $\JML^{k}$.
However, since the surrogate might only be reliable in a specific part of the set of feasible control vectors around the current iterate $\bu_k$,
we retrain the surrogate if the FOM optimization step suggests that a further improvement of the objective function value is possible.
Therefore, the overall goal of the \ROMEnOpt\ algorithm is to terminate with a considerably smaller number of (outer)
iterations $k$ than the \FOMEnOpt\ algorithm, and thus, to reduce the computational costs for solving the polymer
EOR optimization problem in~\eqref{equ:opt_problem}.
We refer to the subsequent section for an extensive complexity and run time comparison for a practical example.
\par
The main motivation for the \ROMEnOpt\ algorithm is
illustrated in \cref{fig:visualization-ROMEnOpt}. Computing the gradient
information using evaluations of the function $\J$ is costly, whereas
gradient computations using the approximation $\JML^k$, obtained,
for instance, via training a neural network, is cheap. In the example, the
\ROMEnOpt\ algorithm performs more optimization steps in total. However,
most of these optimization steps are cheap since they only require
evaluations of $\JML^k$. For the \ROMEnOpt\ algorithm, only those
steps involving evaluations of $\J$ (i.e., outer iterations) require a
large computational effort. Each optimization step is costly in the \FOMEnOpt\ algorithm
since the exact objective function~$\J$ is evaluated
multiple times. Altogether, in the example shown in
\cref{fig:visualization-ROMEnOpt}, the \ROMEnOpt\ algorithm performs less
costly gradient computations than the \FOMEnOpt\ procedure while arriving approximately
at the same optimum. This motivates why the \ROMEnOpt\ algorithm can be
preferable with respect to the required computation time. 

\def\shift{13.5}

\begin{figure}[t]
    \centering
    \resizebox{0.9\textwidth}{!}{
        \begin{tikzpicture}[fom/.style={draw, circle, fill, inner sep=3pt},
                            fom-train/.style={draw, circle, fill, inner sep=4pt, outer sep=2pt, double distance=2pt},
                            rom/.style={draw, circle, inner sep=2pt},
                            fom-edge/.style={-latex, draw, thick},
                            rom-edge/.style={-latex, draw, thick, dotted},
                            fom-rom-edge/.style={-latex, draw, thick, dashed}]
            \node[fom] (start) at (0., 0.) {};
            \node[fom] (fom1) at (3., -2.) {};
            \node[fom] (fom2) at (7., -1.) {};
            \node[fom] (fom3) at (9., -3.5) {};
            \node[fom] (fom4) at (10., -6.) {};
            \node[fom] (fom5) at (10.5, -9.) {};
            \node[fom] (end) at (10., -12.) {};

            \path[fom-edge] (start)
                edge (fom1) (fom1)
                edge (fom2) (fom2)
                edge (fom3) (fom3)
                edge (fom4) (fom4)
                edge (fom5) (fom5)
                edge (end);

            \draw[-, thick] (12.,0.) -- (12.,-12.);

            \node[fom-train] (start-rom) at ($(start)+(\shift,0.)$) {};
            \node[rom] (rom1) at ($(fom1)+(\shift,0.)$) {};
            \node[rom] (rom2) at ($(5.,-1.75)+(\shift,0.)$) {};
            \node[rom] (rom3) at ($(6.5,-3.)+(\shift,0.)$) {};
            \node[rom] (rom4) at ($(8.5,-3.25)+(\shift,0.)$) {};
            \node[fom-train] (fom-rom1) at ($(9.5,-5.)+(\shift,0.)$) {};
            \node[rom] (rom5) at ($(9.,-7.5)+(\shift,0.)$) {};
            \node[rom] (rom6) at ($(8.5,-9.5)+(\shift,0.)$) {};
            \node[rom] (rom7) at ($(9.5,-10.5)+(\shift,0.)$) {};
            \node[fom] (end-rom) at ($(end)+(\shift,0.)$) {};

            \path[rom-edge] (start-rom)
                edge[fom-rom-edge] (rom1) (rom1)
                edge (rom2) (rom2)
                edge (rom3) (rom3)
                edge (rom4) (rom4)
                edge (fom-rom1) (fom-rom1)
                edge[fom-rom-edge] (rom5) (rom5)
                edge (rom6) (rom6)
                edge (rom7) (rom7)
                edge (end-rom);

            \matrix[draw, above right, scale=0.5] at ($(current bounding box.south west)+(-1.,1.)$) {
                \node[fom] at (0.24, 0.275) {};
                \node at (0.75, 0.) {\Large Gradient computation using FOM}; \\
                \node[fom-train] at (0.15, 0.45) {};
                \node[align=left] at (0.75, 0.) {\Large Gradient computation using FOM\\\Large and neural network training}; \\
                \node[rom] at (0.25, 0.325) {};
                \node at (0.75, 0.) {\Large Gradient computation using surrogate}; \\
                \draw (0., 0.375) edge[fom-edge] (0.75, 0.375);
                \node at (0.75, 0.) {\Large FOM optimization step}; \\
                \draw (0., 0.55) edge[fom-rom-edge] (0.75, 0.55);
                \node[align=left] at (0.75, 0.) {\Large Surrogate-based optimization step\\\Large similar to FOM optimization step}; \\
                \draw (0., 0.375) edge[rom-edge] (0.75, 0.375);
                \node at (0.75, 0.) {\Large Surrogate-based optimization step}; \\
            };
        \end{tikzpicture}
    }
	\caption{Example of optimization paths taken by the \FOMEnOpt\ algorithm (left part of the figure) and the \ROMEnOpt\ algorithm (right part of the figure).}
	\label{fig:visualization-ROMEnOpt}
\end{figure}
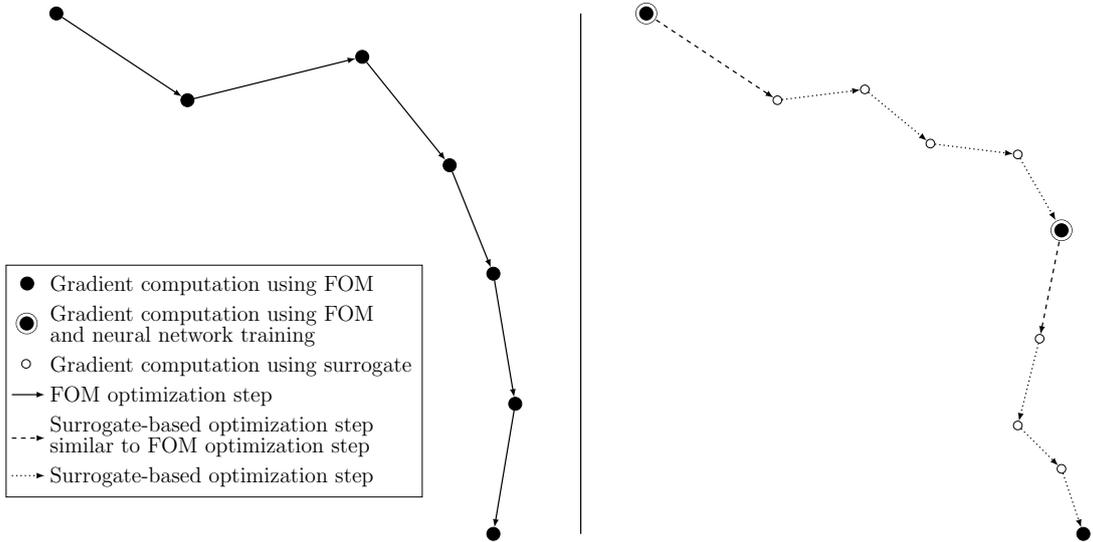

\section{Numerical validation for a five-spot benchmark problem}\label{sec:numerical_experiments}\label{sec:6}
In this section, we present an example with a synthetic oil reservoir in which the polymer flooding optimization
problem~\eqref{equ:opt_problem} is solved using the traditional solution method, the \FOMEnOpt\ algorithm,
and our proposed \ROMEnOpt\ method presented in \cref{alg:ROM_EnOpt}.
The focus is to demonstrate a more efficient and improved method of dealing with the
optimization part of a closed-loop reservoir workflow \cite{chen2009efficient}
for polymer flooding with the assumption that the geological properties of the reservoir are known.
We start by providing information on the algorithm implementation.

\subsection{Implementational details}\label{sec:implementational-details}

For a numerical approximation of the system~\eqref{equ:polymer_model} of non-linear partial differential equations and the corresponding well equations~\eqref{equ:well_equation}, we make use of the open porous media flow reservoir simulator (OPM)~\cite{flo2019open,opm2021}. 
The system is discretized spatially using a two-point flux approximation (TPFA) with upstream-mobility weighting (UMW) and temporally using a fully-implicit Runge-Kutta method.
The resulting discrete-in-time equations are solved using
a Newton-Raphson scheme to obtain time-dependent states
and the output quantities from the well's equation in terms
of fluid production of the reservoir per time step.
In this numerical experiment, we perform all
polymer flooding simulations in parallel on a 50 core CPU.
\par
For the implementation of the DNN-based surrogates, the Python package pyMOR~\cite{milk2016} is used.
The implementation of the neural networks and corresponding training
algorithms in pyMOR is based on the machine learning library
PyTorch~\cite{paszke2019}.
\par
Throughout our numerical experiments described in the subsequent section,
we apply the L-BFGS optimizer with strong Wolfe line-search~\cite{wolfe1969,wolfe1971} for training
the neural networks, i.e., to solve~\eqref{equ:DNNOptimizationProblem}.
Further, we perform a maximum of $1000$ training epochs in each restart.
\par
The number of training restarts influences the accuracy of the trained
neural networks and the computation time required for the training.
A larger number of restarts typically leads to smaller losses and more
training time. To take these two factors into account, we consider
different numbers of restarts in our numerical study presented below.
The respective results can be found in the subsequent section.
In general, we use relatively small numbers of restarts. First of all,
we are not interested in obtaining a neural network with very high
accuracy. Due to the adaptive retraining of the networks, the surrogates
are replaced in each outer iteration anyway. They are only supposed to
lead the optimizer to a point with a larger objective function value.
On the other hand, as indicated before, a larger number of restarts
might result in an unnecessarily long training phase, which must be performed
in each outer iteration. The small numbers of $15$ and $35$ restarts we tried in our
studies can thus be seen as a compromise between the accuracy of the surrogate models
and computational effort for the training algorithm.
\par
We use $10$\% of the sample set for validation during the neural network
training, and the training routine is stopped early if the loss does not
decrease for $10$ consecutive epochs. Moreover, the mean squared
error loss (MSE loss) is used as the loss function.
The neural network training is performed on scaled data.
The input values are scaled to $[0,1]^{N_u}$, and the output values are scaled to $[0,1]$ in the \ANNs-case and $[0,1]^{N_t}$ in the \ANNv-case, respectively.
The scaling of the input values can be computed exactly using the lower and
upper bounds $u_j^\text{low}$ and $u_j^\text{upp}$ for the control variables by
\begin{align}\label{equ:scaling}
	u_j^i \mapsto \frac{u_j^i-u_j^\text{low}}{u_j^\text{upp}-u_j^\text{low}}
\end{align}
for $j=1,\dots,N_w$ and $i=1,\dots,N_t$. For the output values, we take
the minimum and maximum value over the training set as lower and upper bound
and perform the same scaling as in \cref{equ:scaling}.
The $\tanh$ function serves as the activation function for each layer.
Kaiming initialization is applied for initializing the neural network weights.
\par
The input and output dimensions of the neural networks were
already described in \cref{sec:ann_enopt_algorithm} and are different for
the \ANNs- and \ANNv-case.
Regarding the training data for the vector valued case \ANNv, we note that
we require $T_k$ to store $\{(\bu_{k,j},j(\bu_{k,j}))\}_{j=1}^N$ instead
of $\{(\bu_{k,j},\J(\bu_{k,j}))\}_{j=1}^N$, which we did not include in
\cref{alg:OptStep} for brevity.

\subsection{Case study: five-spot field}
\begin{figure}[h]
	\centering
	\includegraphics[width=0.55\textwidth]{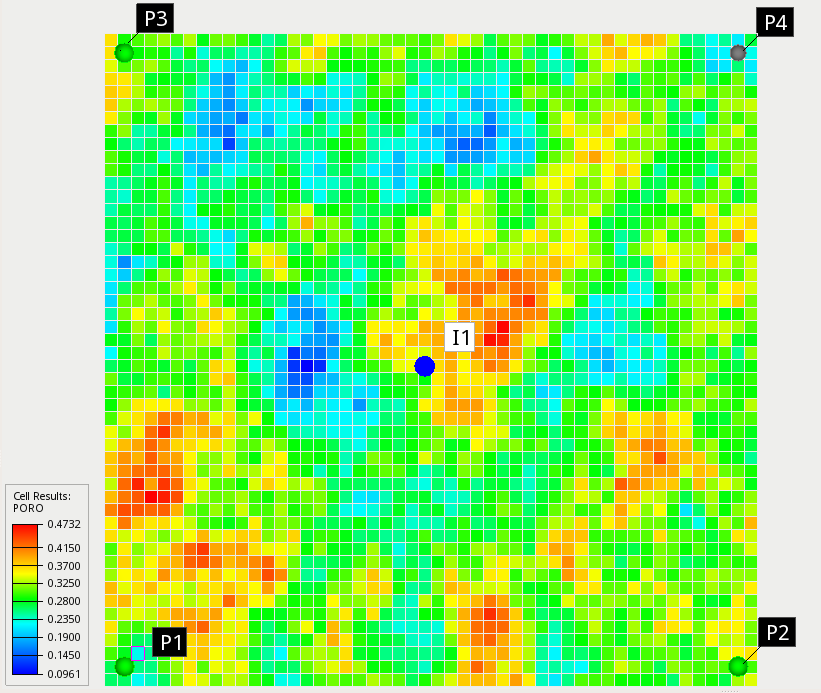}
	\caption{Porosity distribution of the five-spot field and placement of the injection and production wells.}%
	\label{fig:5spotporo}
\end{figure}

The numerical experiment considers a two-dimensional reservoir model with a three-phase flow, including oil, water, and gas (cf. Section \ref{sec:2}).
The computations are performed on a uniform grid that consists of $50\times 50$ grid cells.
The model has one injection and four production wells spatially arranged in a
five-spot pattern as shown in \cref{fig:5spotporo}.

On average, the reservoir has approximately $30$\% porosity with a
heterogeneous permeability distribution. The initial reservoir pressure is $200$~bar. The initial average oil
and water saturations are $0.6546$ and $0.3454$, respectively. The original oil in place is
$4.983\cdot 10^6$~sm$^3$. Fluid properties are similar to those of a light oil reservoir.
The viscosity for saturated oil at varying bubble point pressure lies
between $0.1$~cP and $0.56$~cP, and the viscosity of water is $0.01$~cP.
The densities of oil and water are taken as $732$~kg/m$^3$ and $1000$~kg/m$^3$, respectively.
In this setting, it is easy to see that the displacement is unfavorable since the oil-water mobility ratio $\lambda$
is such that $10\leq\lambda\leq 56$. The reservoir rock parameters utilized for the polymer flooding simulation in this
problem are given by \cref{tab:Rock1}.

\begin{table}[h]
	\centering
	\begin{tabular}{l c c}
		\toprule
		Parameter &  {Value}  & {Unit} \\ \midrule \midrule
		Dead pore space for polymer solution  & {$0.1800$} & {$-$}   \\  
		Maximum polymer adsorption value & {$7.5\cdot 10^{-4}$} & {kg/kg}   \\  
		Residual resistance factor of polymer solution & {$2.5$} &  {$-$}  \\ 
		Reservoir rock density & {$1980$} &  {$\text{kg/rm}^3$}  \\ 
		Polymer mixing parameters  & {$0.65$} &  {$-$}  \\ \bottomrule
	\end{tabular}
	\caption{Reservoir model parameters used in the polymer flooding simulations.}
	\label{tab:Rock1}
\end{table}

In this example, the injection well is controlled by two independent control
variables, namely the water injection rate and the polymer concentration at each control
time step. The lower and upper bounds for the water injection rate are set to $0$~sm$^3$/day
and $2000$~sm$^3$/day respectively, while the lower and upper bounds for the polymer concentration are
set to $0$~kg/sm$^3$ and $2.5$~kg/sm$^3$. Hence, the polymer injection rate
ranges from $0$ to $5000$~kg/day. Each production well is controlled by a reservoir
fluid production rate target with a lower limit of $0$~sm$^3$/day and an upper limit of $500$~sm$^3$/day.
Bottom hole pressure limits are imposed on the wells, namely a maximum of $500$~bar for the injector and a
minimum of $150$~bar for each producer. The production period for the reservoir is set to $50$~months, and the
control time step is taken as $5$~months. Therefore, there are $N_u = (2+4)\times 10 = 60$ control variables in total
to solve for in~\eqref{equ:opt_problem}. For the objective function~\eqref{equ:fom-objective},
we used the economic parameters listed in \cref{tab:Economy1}.

\begin{table}[h]
	\centering
	\begin{tabular}{l c c}
		\toprule
		Parameter &  {Value}  & {Unit} \\ \midrule \midrule
		Oil price $r_{\text{OP}}$ & {$500$} & {USD/sm$^3$} \\
		Price of gas production $r_{\text{GP}}$ & {$0.15$} & {USD/sm$^3$} \\
		Cost of polymer injection $r_{\text{PI}}$  & {2.5} & {USD/kg} \\
		Cost of polymer production $r_{\text{PP}}$  & {0.5} & {USD/kg} \\
		Cost of  water injection or production $r_{\text{WI}}$, $r_{\text{WP}}$ & {$30$} & {USD/sm$^3$} \\
		Annual discount rate $d_{\tau}$ & {$0.1$} &  {$-$} \\
		\bottomrule
	\end{tabular}
	\caption{Economic parameters used in the numerical experiments.} \label{tab:Economy1}
\end{table}

Using the two different surrogate models for the objective
function~\eqref{equ:fom-objective} constructed by means of neural networks,
namely \ANNs\ and \ANNv\ as explained in \cref{sec:ann_enopt_algorithm}, the optimization
problem~\eqref{equ:opt_problem} is solved using the \ROMEnOpt\ algorithm.
In this case, the \ROMEnOpt\ algorithm for~\eqref{equ:opt_problem} using \ANNs\ and \ANNv\ to
approximate the objective function $\J$ from~\eqref{equ:fom-objective}
is denoted by \ROMEnOpts\ and \ROMEnOptv, respectively.
The \EnOpt\ parameters for both, the \FOMEnOpt\ and the two variants of the
\ROMEnOpt\ method, are presented in \cref{tab:enOpt-para}. We remark that
the tolerances $\varepsilon$, $\varepsilon_i$, and $\varepsilon_o$ are
applied to the scaled quantities, i.e., the output quantities, for which
the respective stopping criteria in \cref{alg:EnOpt,alg:ROM_EnOpt} are
checked, have already been scaled as described in
\cref{sec:implementational-details}.

\begin{table}[h]
	\centering
	\begin{tabular}{l c}
		\toprule
		Parameter &  {Value} \\ \midrule \midrule
		Initial step size $\beta_0$ & {$0.3$} \\
		Step size contraction $r$ & {$0.5$} \\
		Maximum step size trials $\nu^*$ & {$10$} \\
		Initial control-type variance $\sigma_{j}$ & {$0.001$} \\
		Constant correlation factor $\rho$ & {$0.9$} \\
		Perturbation size $N$ & {$100$} \\
		\midrule
		Tolerances \begin{tabular}{c}
			\FOMEnOpt\ $\varepsilon$ \\
			\ROMEnOpt\ inner iteration $\varepsilon_i$ \\
			\ROMEnOpt\ outer iteration $\varepsilon_o$
		\end{tabular} &
		\begin{tabular}{c}
			$10^{-6}$ \\
			$10^{-6}$ \\
			$10^{-2}$
		\end{tabular} \\
		\bottomrule
	\end{tabular}
	\caption{Parameters used in the \FOMEnOpt\ and \ROMEnOpt\ algorithms.}
	\label{tab:enOpt-para}
\end{table}

We compare the \ROMEnOpt\ results with those of the \FOMEnOpt\ algorithm for two different initial
guesses $\bu_0^1\in\mathcal{D}_\text{ad}$ and $\bu_0^2\in\mathcal{D}_\text{ad}$.
The initial solution $\bu_0^1$ includes $700$~sm$^3$/day
for the water injection rate at the injection well, $150$~sm$^3$/day for the reservoir fluid production
rate at each production well, and $0.5$~kg/sm$^3$ for the polymer concentration
(equivalently $350$~kg/day for polymer injection rate) at the injection well over the simulation period.
Similarly, $\bu_0^2$ includes $600$~sm$^3$/day for the water injection rate, $100$~sm$^3$/day for the reservoir
fluid production rate, and $0.5$~kg/sm$^3$ for the polymer concentration.

\begin{figure}[h]
	\centering
	\input{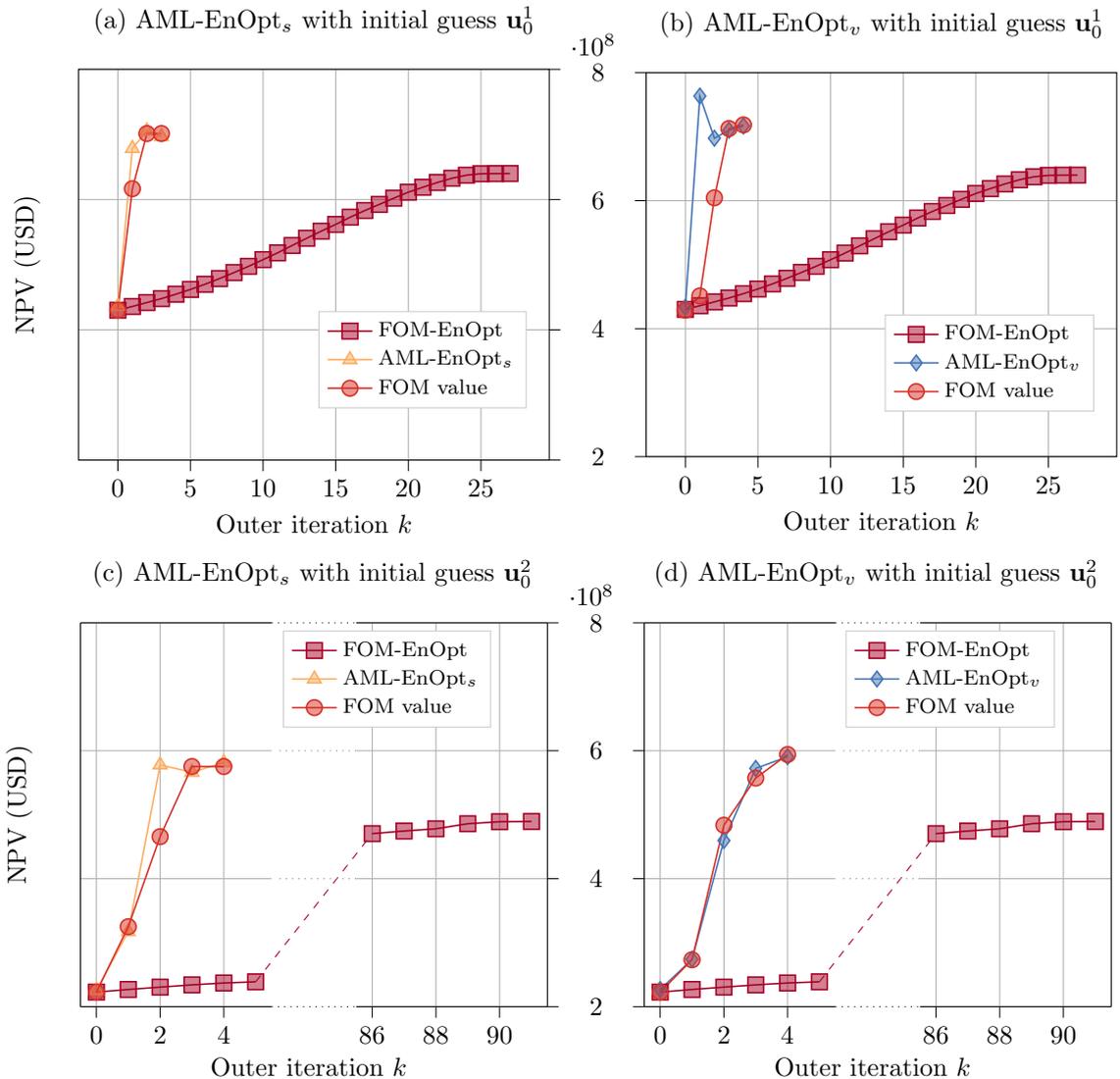}
	\caption{Comparison of the NPV values obtained during the outer iterations of the \FOMEnOpt, \ROMEnOpts, and \ROMEnOptv\ procedures for two different initial guesses $\bu_0^1\in\mathcal{D}_\text{ad}$ and $\bu_0^2\in\mathcal{D}_\text{ad}$. For each \ROMEnOpt\ procedure, the corresponding FOM value $\J(\bu_k)$ at the current iterate $\bu_k$ of the respective \ROMEnOpt\ method is indicated as well.}
	\label{fig:5Spot-diff-ini-guess}
\end{figure}

\cref{fig:5Spot-diff-ini-guess} compares the values of the objective
function during the outer iterations of the
\FOMEnOpt, \ROMEnOpts, and \ROMEnOptv\ strategies using the initial
solutions~$\bu_0^1$ and~$\bu_0^2$. Furthermore, the value~$\J(\bu_k)$
at the outer iterate $\bu_k$ (denoted by ``FOM value'') for the
respective \ROMEnOpt\ method is depicted.
\par
Since the \ROMEnOpt\ algorithms only use an approximate surrogate
model $\JML^k$, the values of $\J$ and $\JML^k$ are not
necessarily the same for the control $\bu_k$. This behavior is especially apparent in
\cref{fig:5Spot-diff-ini-guess}(b), where the \ROMEnOptv\ algorithm is
examined for the initial guess $\bu_0^1$. Here, after the first outer
iteration, the values $\J(\bu_1)$ and $\JML^0(\bu_1)$ differ from
each other by a significant amount. A possible reason is that the surrogate model~$\JML^0$
does not extrapolate well to the region where the first (inner)
\ROMEnOpt\ iteration converged to. This further indicates that the found iterate~$\bu_1$ is far 
from the initial solution $\bu_0$, where the initial model~$\JML^0$ was trained.
However, since the \ROMEnOpt\ algorithm uses
evaluations of $\J$ in the stopping criterion, the \ROMEnOpt\ does not terminate but
continues by training a new surrogate model using training data sampled normally
around~$\bu_1$. Hence, the new surrogate $\JML^1$ tries to approximate the
objective function $\J$ well around $\bu_1$. In each plot, we see that in the
last two outer iterations of the respective \ROMEnOpt\ procedure, the FOM
value and the \ROMEnOpt\ value agree to minimal deviations. This suggests
that the surrogate model approximates the full objective function
well in the region of the (local) optimum found by the \ROMEnOpt\ method.
\par
More so, in \cref{fig:5Spot-diff-ini-guess}, it is seen that both, the \ROMEnOpts\
and the \ROMEnOptv\ algorithm, require considerably less (costly) outer
iterations than the \FOMEnOpt\ method. This leads to an improvement in the
run time of the method, which is detailed in \cref{tab:resultcomp1}.
Besides the faster convergence of the method, we also remark that the
\ROMEnOpt\ algorithms find local optima with larger
objective function values than the \FOMEnOpt\ algorithm. However, since the
objective function $\J$ is multi-modal, this is not guaranteed.
\par
We emphasize that each outer iteration of the \ROMEnOpt\ algorithm includes
many inner iterations (see also \cref{tab:resultcomp1,tab:resultcomp2}),
which leads to the large jumps in the objective function values between
consecutive outer iterations, as present in \cref{fig:5Spot-diff-ini-guess}.

Further comparisons in terms of function values, numbers of inner and
outer iterations, numbers of evaluations of the FOM function $\J$ and
surrogate approximations $\JML^k$, total run time, and speedup are
presented in \cref{tab:resultcomp1,tab:resultcomp2}.
\par
With the different initial guesses $\bu_0^1$ and $\bu_0^2,$ we found that the
number of outer iterations required by the \FOMEnOpt\ algorithm significantly
differs. However, the \ROMEnOpt\ methods require only~$4$ and~$5$ outer
iterations. This reduced number of outer iterations leads to a remarkable
speedup in the overall computation time $T_\text{total}$ and is particularly
reflected in the reduced number of FOM evaluations, i.e., evaluations of the
objective function~$\J$, which require costly polymer flooding simulations.
Although each outer iteration consists of multiple inner
iterations using the surrogate $\JML^k$, it does not
contribute substantially to the overall run time because evaluating the
surrogate $\JML^k$ is very cheap.

\begin{table}[h]
	\centering
	\resizebox{\columnwidth}{!}{%
		\begin{tabular}{l c c c c c c c c c c c}
			\toprule
			{Method} & {FOM value} & {Surrogate value} & {Outer it.} & {Inner it.}& {FOM eval.} & {Surrogate eval.} & {$T_\text{total}$ (min)} & {Speedup} \\ \midrule \midrule
			
			\FOMEnOpt & {$6.400\cdot 10^8$} & {$-$} & {28}&{$-$} & {$2839$} & {$-$}& {$54.86$}& {$-$}\\ 
			
			\ROMEnOpts & {$7.013\cdot 10^8$} & {$6.968\cdot 10^8$} & {$4$}&{$233$} & {$407$} & {$12315$}& {$8.87$}& {$6.18$}\\ 
			
			\ROMEnOptv & {$7.185\cdot 10^8$} & {$7.168\cdot 10^8$} & {$5$}& {$312$}&{$509$} & {$14101$} & {$14.10$}& {$3.89$}\\
			\bottomrule
	\end{tabular}}
	\caption{Comparisons of the results from the different solution strategies \FOMEnOpt, \ROMEnOpts, and \ROMEnOptv\ using the initial guess $\bu_0^1$ and $N_1=N_2=35$ neurons in each hidden layer and $15$ restarts for the neural network training.}
	\label{tab:resultcomp1}
\end{table}

\begin{table}[!h]
	\centering
    \resizebox{\columnwidth}{!}{%
    	\begin{tabular}{l c c c c c c c c c c c}
    		\toprule
    		{Method} & {FOM value} & {Surrogate value} & {Outer it.} & {Inner it.}& {FOM eval.} & {Surrogate eval.} & {$T_\text{total}$ (min)} & {Speedup} \\ \midrule \midrule
		
	    	\FOMEnOpt & {$4.895\cdot 10^8$} & {$-$} & {92}&{$-$} & {$9310$} & {$-$}& {$135.77$}& {$-$}\\ 
		
		    \ROMEnOpts & {$5.754\cdot 10^8$} & {$5.816\cdot 10^8$} & {$4$}&{$111$} & {$407$} & {$10837$}& {$9.85$}& {$13.78$}\\ 
		
    		\ROMEnOptv & {$5.942\cdot 10^8$} & {$5.908\cdot 10^8$} & {$4$}& {$117$}&{$407$} & {$10033$} &{$11.05$}& {$12.29$}\\
	    	\bottomrule
    \end{tabular}}
	\caption{Comparisons of the results from the different solution strategies \FOMEnOpt, \ROMEnOpts, and \ROMEnOptv\ using the initial guess $\bu_0^2$ and $N_1=N_2=25$ neurons in each hidden layer and $35$ restarts for the neural network training.}
    \label{tab:resultcomp2}
\end{table}

\begin{figure}[h!]
	\centering
	\input{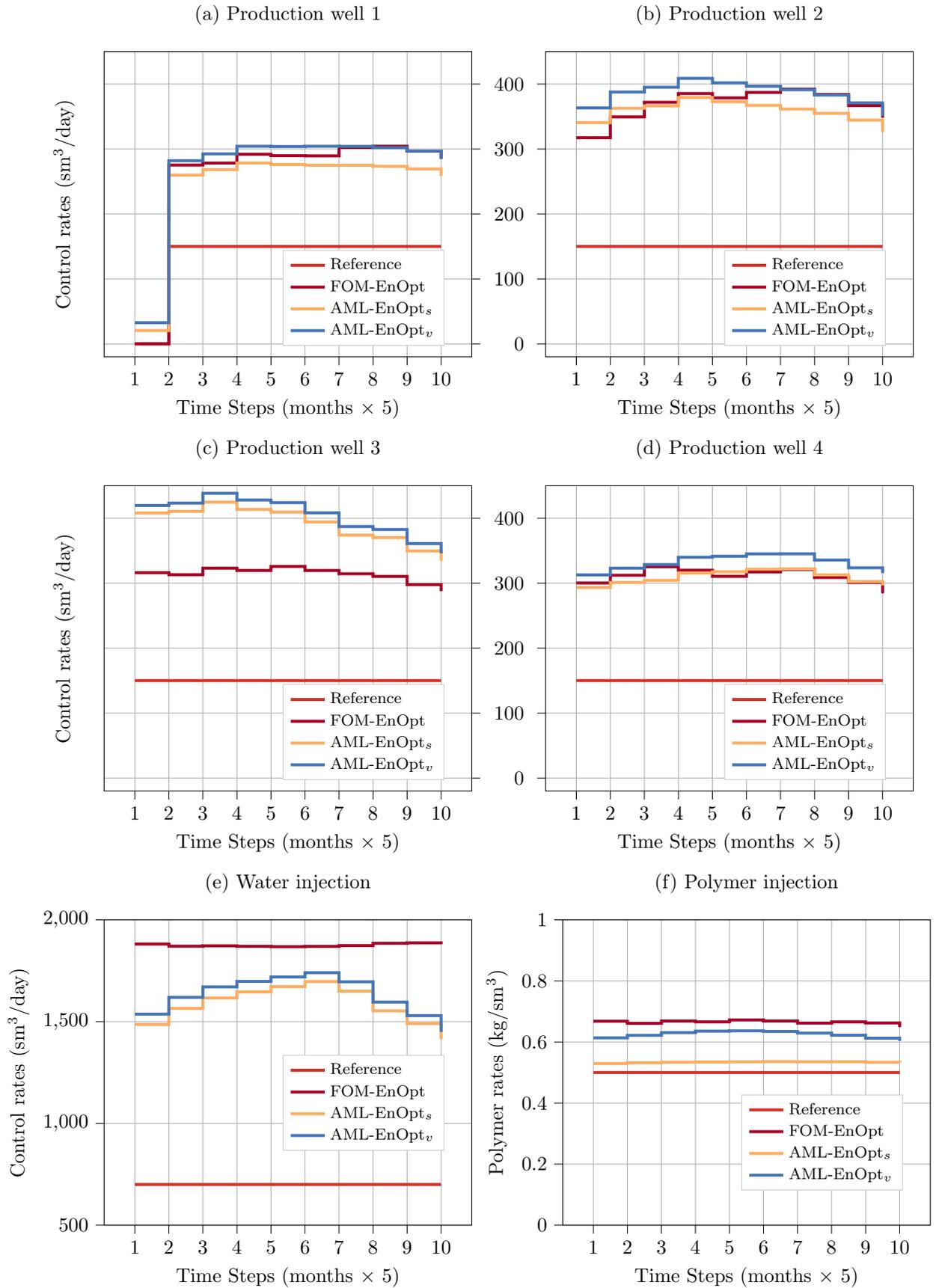}
	\caption{Comparison of the optimal solutions obtained via the \FOMEnOpt, \ROMEnOpts, and \ROMEnOptv\ algorithms using the initial guess $\bu_0^1$, which is depicted as the reference solution.}
	\label{fig:5Spot-control}
\end{figure}

For the initial solution $\bu_0^1$, the optimizers obtained from the
three solution strategies are depicted in \cref{fig:5Spot-control}.
Further, the initial guess $\bu_0^1$ is shown as a reference.
\par
The control variables obtained by the \ROMEnOpts\ and the \ROMEnOptv\
algorithm are close to those of the \FOMEnOpt\ method, except for production
well 3 (see \cref{fig:5Spot-control}(c)) and the water injection rate (see
\cref{fig:5Spot-control}(e)). For each control variable, the values obtained
via the \ROMEnOpts\ and \ROMEnOptv\ procedures are close to each other.
Together with the FOM values of \ROMEnOpts\ and \ROMEnOptv\ presented in
\cref{tab:resultcomp1} and the evolution of the FOM values for the two
methods shown in \cref{fig:5Spot-diff-ini-guess}(a)-(b), this suggests that
the \ROMEnOpts\ and the \ROMEnOptv\ methods traverse almost the same path in
the control space $\mathcal{D}_\text{ad}$ and find local optima close to each other.

\begin{figure}[h!]
	\centering
	\begin{tikzpicture}
\begin{axis}[%
name=left,
anchor=west,
width=10cm,
scale=0.75,
legend cell align={left},
legend style={font=\footnotesize, fill opacity=1, draw opacity=1, text opacity=1, xshift=0pt, yshift=-105pt, draw=white!80!black},
tick align=outside,
tick pos=left,
x grid style={white!69.0196078431373!black},
xlabel={Times (Months $\times$ 5)},
xmajorgrids,
xtick style={color=black},
xtick distance=1,
y grid style={white!69.0196078431373!black},
ymajorgrids,
ymin=20000, ymax=1600000,
ylabel={FOPT (sm$^3$)},
ytick pos=left,
yticklabel pos=left,
ytick style={color=black},
tick scale binop=\cdot
]
\addplot[semithick, \FOMEnOptColor, mark=\FOMEnOptMarker, mark size=3, mark options={solid, fill opacity=0.5}] table[row sep=crcr]{%
x	y\\
1	28946.423828125\\
2	178810.125\\
3	371541.3125\\
4	568178.4375\\
5	769268.375\\
6	947956\\
7	1079331.375\\
8	1172276\\
9	1242734.75\\
10	1342491.875\\
};
\addlegendentry{\FOMEnOpt}

\addplot[semithick, \ROMEnOptsColor, mark=\ROMEnOptsMarker, mark size=3, mark options={solid, fill opacity=0.5}] table[row sep=crcr]{%
x	y\\
1	31271.75390625\\
2	187528\\
3	359340.75\\
4	523606.5625\\
5	691319.6875\\
6	857118.1875\\
7	1018350\\
8	1162600.5\\
9	1277819.5\\
10	1425193.75\\
};
\addlegendentry{\ROMEnOpts}

\addplot [semithick, \ROMEnOptvColor, mark=\ROMEnOptvMarker, mark size=3, mark options={solid, fill opacity=0.5}]
table[row sep=crcr]{%
	1	34970.16796875\\
	2	213107.015625\\
	3	427347.90625\\
	4	633028.25\\
	5	841052.5625\\
	6	1039965.625\\
	7	1215954.5\\
	8	1349180.375\\
	9	1440156\\
	10	1554116.375\\
};
\addlegendentry{\ROMEnOptv}

\end{axis}

\begin{axis}[%
name=right,
anchor=west,
width=10cm,
at=(left.east),
xshift=2cm,
scale=0.75,
legend cell align={left},
legend style={font=\footnotesize, fill opacity=1, draw opacity=1, text opacity=1, xshift=-60pt, draw=white!80!black},
tick align=outside,
tick pos=left,
x grid style={white!69.0196078431373!black},
xlabel={Times (Months $\times$ 5)},
xmajorgrids,
xtick style={color=black},
xtick distance=1,
y grid style={white!69.0196078431373!black},
ymajorgrids,
ymin=0, ymax=600000,
yticklabel pos=left,
ytick style={color=black},
ylabel={FWPT (sm$^3$)},
tick scale binop=\cdot
]
\addplot[semithick, \FOMEnOptColor, mark=\FOMEnOptMarker, mark size=3, mark options={solid, fill opacity=0.5}] table[row sep=crcr]{%
x	y\\
1	0.603970646858215\\
2	0.603970646858215\\
3	0.603970646858215\\
4	0.603970646858215\\
5	0.603970646858215\\
6	19883.671875\\
7	89967.453125\\
8	199798.765625\\
9	325304.46875\\
10	561535\\
};
\addlegendentry{\FOMEnOpt}

\addplot[semithick, \ROMEnOptsColor, mark=\ROMEnOptsMarker, mark size=3, mark options={solid, fill opacity=0.5}] table[row sep=crcr]{%
x	y\\
1	1.54989981651306\\
2	8.87867641448975\\
3	19.2297782897949\\
4	29.7683181762695\\
5	40.5210227966309\\
6	96.8332824707031\\
7	2050.07958984375\\
8	17060.392578125\\
9	63279.12109375\\
10	213726.890625\\
};
\addlegendentry{\ROMEnOpts}

\addplot[semithick, \ROMEnOptvColor, mark=\ROMEnOptvMarker, mark size=3, mark options={solid, fill opacity=0.5}] table[row sep=crcr]{%
	x	y\\
1	1.55533587932587\\
2	5.39974212646484\\
3	13.1384220123291\\
4	22.4992027282715\\
5	32.784912109375\\
6	6664.89697265625\\
7	39470.29296875\\
8	119321.9375\\
9	238619.078125\\
10	485332.0625\\
};
\addlegendentry{\ROMEnOptv}

\end{axis}
\node[anchor=south, yshift=4pt] at (left.north) {(a) Oil production};
\node[anchor=south, yshift=4pt] at (right.north) {(b) Water production};

\begin{axis}[%
name=bottom,
anchor=north,
width=10cm,
at=(left.south),
yshift=-59pt,
xshift=.25\textwidth,
scale=0.75,
legend cell align={left},
legend style={font=\footnotesize, fill opacity=1, draw opacity=1, text opacity=1, xshift=-60pt, draw=white!80!black},
tick align=outside,
tick pos=left,
x grid style={white!69.0196078431373!black},
xlabel={Times (Months $\times$ 5)},
xmajorgrids,
xtick style={color=black},
xtick distance=1,
y grid style={white!69.0196078431373!black},
ymajorgrids,
ymin=0, ymax=5100,
ylabel={FCPT (kg)},
ytick pos=left,
yticklabel pos=left,
ytick style={color=black}
]
\addplot[semithick, \FOMEnOptColor, mark=\FOMEnOptMarker, mark size=3, mark options={solid, fill opacity=0.5}] table[row sep=crcr]{%
	x	y\\
1	0\\
2	0\\
3	0\\
4	0\\
5	0\\
6	9.37908172607422\\
7	100.712509155273\\
8	457.565216064453\\
9	1335.78076171875\\
10	5010.30419921875\\
};
\addlegendentry{\FOMEnOpt}

\addplot[semithick, \ROMEnOptsColor, mark=\ROMEnOptsMarker, mark size=3, mark options={solid, fill opacity=0.5}] table[row sep=crcr]{%
	x	y\\
	1	0\\
	2	0\\
	3	0\\
	4	0\\
	5	1.21638961948378e-26\\
	6	4.58699616956437e-07\\
	7	0.130956009030342\\
	8	4.55327558517456\\
	9	40.8345489501953\\
	10	416.389251708984\\
};
\addlegendentry{\ROMEnOpts}

\addplot[semithick, \ROMEnOptvColor, mark=\ROMEnOptvMarker, mark size=3, mark options={solid, fill opacity=0.5}] table[row sep=crcr]{%
	x	y\\
	1	0\\
	2	0\\
	3	0\\
	4	0\\
	5	3.22611446463537e-25\\
	6	3.51576161384583\\
	7	24.9536437988281\\
	8	125.765975952148\\
	9	436.031707763672\\
	10	2201.5537109375\\
};
\addlegendentry{\ROMEnOptv}

	\end{axis}
\node[anchor=south, yshift=4pt] at (bottom.north) {(c) Polymer production};
\end{tikzpicture}%
	\caption{Comparison of the production data obtained from the different solution strategies \FOMEnOpt, \ROMEnOpts, and \ROMEnOptv\ using the initial guess $\bu_0^1$.}
	\label{fig:5Spot-prodData}
\end{figure}
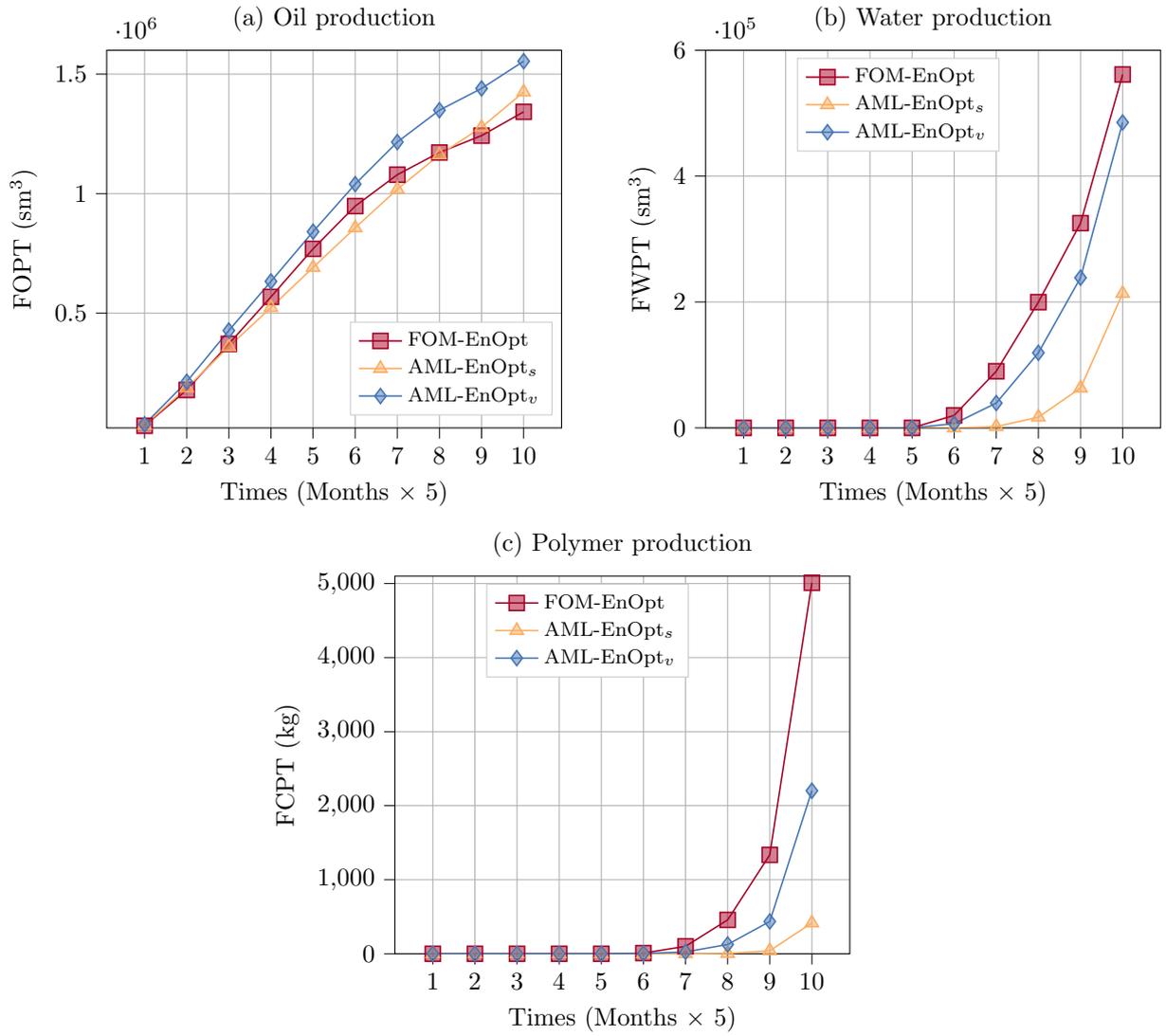

\cref{fig:5Spot-prodData}(a) depicts a comparison of the total field oil production for the optimal solutions
(in \cref{fig:5Spot-control}) of the three solution methods. The total field oil production by \FOMEnOpt,
\ROMEnOpts, and \ROMEnOptv\ are $1.343\cdot 10^6$, $1.425\cdot 10^6$, and $1.554\cdot 10^6$ (in sm$^3$),
respectively. The solution obtained by \ROMEnOptv\ attains the highest oil production in total,
followed by the \ROMEnOpts. The total back-produced water and polymer from operating the five-spot field
with the different optimal solutions are presented in \cref{fig:5Spot-prodData}(b) and \cref{fig:5Spot-prodData}(c), respectively.
Here, we found that the \ROMEnOpts\ and \ROMEnOptv\ solutions are more
economical and environmentally friendly than the one provided by using the
\FOMEnOpt\ method.

\begin{figure}[h!]
	\centering
	\begin{tikzpicture}
\begin{axis}[
name=left,
anchor=west,
width=10cm,
scale=0.75,
legend cell align={left},
legend style={font=\footnotesize, fill opacity=1, draw opacity=1, text opacity=1, xshift=-10pt, yshift=-85pt, draw=white!80!black},
tick align=outside,
tick pos=left,
x grid style={white!69.0196078431373!black},
xlabel={Outer iteration $k$},
xmajorgrids,
xmin=-0.1, xmax=5.1,
xtick style={color=black},
scaled y ticks = false,
y grid style={white!69.0196078431373!black},
ymajorgrids,
ymin=3e+8, ymax=9e+8,
ylabel={NPV (USD)},
yticklabels={,,},
ytick pos=right,
yticklabel pos=left,
ytick style={color=black}
]
\addplot [semithick, color0, mark=triangle*, mark size=3, mark options={solid, fill opacity=0.5}]
table {%
	0	408680000
	1	722760000
	2	666520000
	3	700280000
};
\addlegendentry{$N_1=N_2=20$}
\addplot [semithick, color2, mark=*, mark size=3, mark options={solid, fill opacity=0.5}]
table {%
	0	386640000
	1	650960000
};
\addlegendentry{$N_1=N_2=25$}
\addplot [semithick, color3, mark=square*, mark size=3, mark options={solid, fill opacity=0.5}]
table {%
	0	376680000
	1	553800000
	2	627680000
	3	635040000
	4	645280000
	5	658840000
};
\addlegendentry{$N_1=N_2=30$}
\addplot [semithick, color5, mark=diamond*, mark size=3, mark options={solid, fill opacity=0.5}]
table {%
	0	396040000
	1	678480000
	2	705240000
	3	696760000
};
\addlegendentry{$N_1=N_2=35$}
\end{axis}

\begin{axis}[
name=right,
anchor=west,
width=10cm,
at=(left.east),
xshift=1.0cm,
scale=0.75,
legend cell align={left},
legend style={font=\footnotesize, fill opacity=1, draw opacity=1, text opacity=1, xshift=-10pt, yshift=-85pt, draw=white!80!black},
tick align=outside,
tick pos=left,
x grid style={white!69.0196078431373!black},
xlabel={Outer iteration $k$},
xmajorgrids,
xmin=-0.1, xmax=4.1,
xtick style={color=black},
y grid style={white!69.0196078431373!black},
ymajorgrids,
ymin=3e+8, ymax=9e+8,
ytick pos=left,
ytick style={color=black},
yticklabel shift=4pt,
y tick scale label style={xshift=-28pt},
tick scale binop=\cdot
]
\addplot [semithick, color0, mark=triangle*, mark size=3, mark options={solid, fill opacity=0.5}]
table {%
0	444480000
1	888040000
2	653200000
3	693400000
4	705120000
};
\addlegendentry{$N_1=N_2=20$}
\addplot [semithick, color2, mark=*, mark size=3, mark options={solid, fill opacity=0.5}]
table {%
0	438040000
1	805680000
2	671600000
3	685200000
4	696800000
};
\addlegendentry{$N_1=N_2=25$}
\addplot [semithick, color3, mark=square*, mark size=3, mark options={solid, fill opacity=0.5}]
table {%
0	410640000
1	689960000
2	675240000
3	685640000
4	700800000
};
\addlegendentry{$N_1=N_2=30$}
\addplot [semithick, color5, mark=diamond*, mark size=3, mark options={solid, fill opacity=0.5}]
table {%
0	411960000
1	763560000
2	697600000
3	710080000
4	716840000
};
\addlegendentry{$N_1=N_2=35$}
\end{axis}
\node[anchor=south, yshift=4pt] at (left.north) {(a) \ROMEnOpts};
\node[anchor=south, yshift=4pt] at (right.north) {(b) \ROMEnOptv};
\end{tikzpicture}
	\caption{Comparison of the \ROMEnOpt\ procedures \ROMEnOpts\ and \ROMEnOptv\ for different numbers of neurons in the hidden layers with fixed initial guess $\bu_0^1$.}
	\label{fig:5Spot-Diff-Neu}
\end{figure}
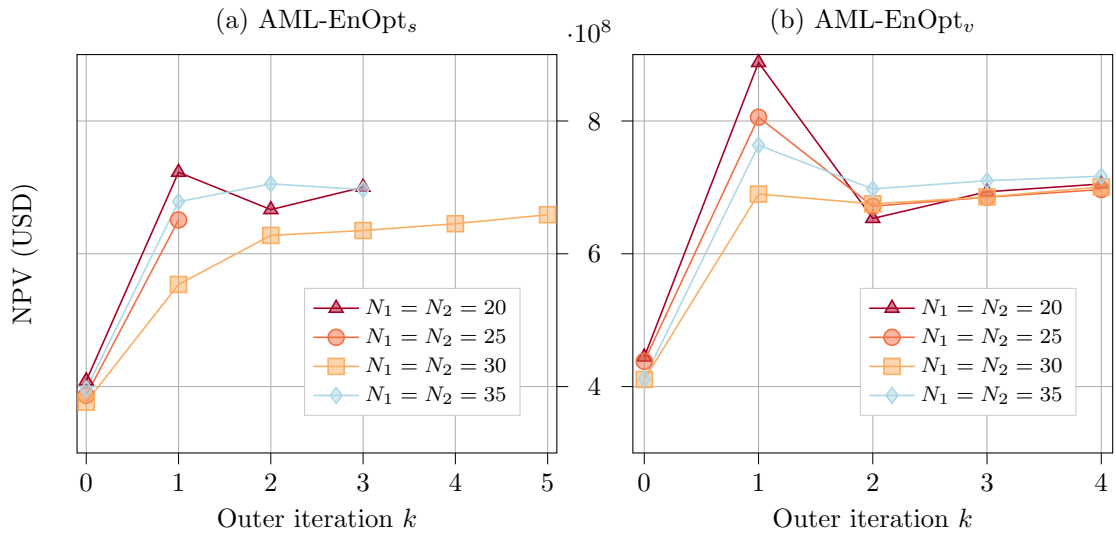

To further investigate the effects of different neural network architectures
on the resulting NPV values, \cref{fig:5Spot-Diff-Neu} depicts the NPV
values obtained by the \ROMEnOpt\ algorithm when using different numbers of
neurons in the hidden layers of the surrogate models \ANNs\ and \ANNv.
\par
We observe that the \ROMEnOptv\ method results are very similar,
which suggests that the \ANNv-approach is more robust and leads
to similar optimal solutions independent of the neural network structure.
In the case of the \ROMEnOpts\ algorithm, different numbers of
neurons lead to results with a larger variation. In particular, the number of
outer iterations performed is different. Hence, the architecture of the
underlying network seems to have a significant effect on the performance of the
resulting \ROMEnOpts\ algorithm.

The maximum, minimum, and average training and validation losses that occurred
in the \ROMEnOpts\ and \ROMEnOptv\ algorithm for the initial guess $\bu_0^1$ are
presented in \cref{tab:losses}. The table shows the respective MSE losses for
different numbers of neurons in the hidden layers.
\par
The results in \cref{tab:losses} do not suggest a significant influence of the number
of neurons on the training and validation results. Further, the scalar- and 
vector-valued cases, \ANNs\ and \ANNv\ respectively, perform similarly in overall training and validation losses. However, we emphasize that, in the
\ANNv\ case, the MSE loss cannot be related directly to the difference in
the output function. Instead, one has to take into account that the outputs of
\ANNv\ are summed up to obtain the surrogate $\JML^k$, while
the MSE loss is measured on the vector-valued outputs of the neural network.
\par
Altogether, the numerical experiments with different numbers of neurons suggest
that already small DNNs with only $20$ neurons in each of the hidden layers yield
appropriate results. In this specific application, we do not benefit from
increasing the complexity of the neural network. We have seen the same behavior
when using more than two hidden layers.
\par
\begin{table}[h!]
	\centering
	\resizebox{\columnwidth}{!}{%
	\begin{tabular}{l c c | c c c | c c c}
		\toprule
		\multirow{2}{*}{Method} & \multirow{2}{*}{\shortstack{Neurons \\ $N_1=N_2$}} & \multirow{2}{*}{\shortstack{Outer\\iter.}} & \multicolumn{3}{c|}{Training loss} & \multicolumn{3}{c}{Validation loss} \\
		& & & Max. & Min. & Avg. & Max. & Min. & Avg. \\ \midrule \midrule
		\ANNs & $20$ & $4$ & $1.2\cdot 10^{-4}$ & $1.3\cdot 10^{-6}$ & $5.3\cdot 10^{-5}$ & $5.4\cdot 10^{-3}$ & $6.7\cdot 10^{-5}$ & $2.3\cdot 10^{-3}$ \\
		\ANNs & $25$ & $2$ & $6.0\cdot 10^{-4}$ & $1.3\cdot 10^{-6}$ & $3.0\cdot 10^{-4}$ & $2.3\cdot 10^{-3}$ & $2.1\cdot 10^{-3}$ & $2.2\cdot 10^{-3}$ \\
		\ANNs & $30$ & $7$ & $7.9\cdot 10^{-4}$ & $8.5\cdot 10^{-7}$ & $1.7\cdot 10^{-4}$ & $6.6\cdot 10^{-3}$ & $1.7\cdot 10^{-3}$ & $3.6\cdot 10^{-3}$ \\
		\ANNs & $35$ & $4$ & $1.8\cdot 10^{-4}$ & $6.1\cdot 10^{-6}$ & $8.2\cdot 10^{-5}$ & $5.5\cdot 10^{-3}$ & $3.7\cdot 10^{-4}$ & $2.7\cdot 10^{-3}$ \\ \midrule
		\ANNv & $20$ & $5$ & $1.8\cdot 10^{-3}$ & $2.1\cdot 10^{-5}$ & $5.2\cdot 10^{-4}$ & $6.9\cdot 10^{-3}$ & $1.2\cdot 10^{-3}$ & $4.2\cdot 10^{-3}$ \\
		\ANNv & $25$ & $6$ & $9.9\cdot 10^{-4}$ & $1.5\cdot 10^{-5}$ & $4.1\cdot 10^{-4}$ & $6.4\cdot 10^{-3}$ & $9.4\cdot 10^{-4}$ & $3.6\cdot 10^{-3}$ \\
		\ANNv & $30$ & $5$ & $9.0\cdot 10^{-4}$ & $9.9\cdot 10^{-6}$ & $4.0\cdot 10^{-4}$ & $1.0\cdot 10^{-2}$ & $5.2\cdot 10^{-4}$ & $4.3\cdot 10^{-3}$ \\
		\ANNv & $35$ & $5$ & $6.0\cdot 10^{-4}$ & $1.1\cdot 10^{-6}$ & $2.2\cdot 10^{-4}$ & $8.8\cdot 10^{-3}$ & $4.6\cdot 10^{-4}$ & $4.3\cdot 10^{-3}$ \\
		\bottomrule
	\end{tabular}}
	\caption{Maximum, minimum, and average MSE loss in the \ROMEnOpts\ and \ROMEnOptv\ algorithm with different numbers of neurons in the hidden layers of the neural networks \ANNs\ and \ANNv\ for fixed initial guess $\bu_0^1$. The number of hidden layers is fixed to two.}
	\label{tab:losses}
\end{table}

\clearpage

\section{Conclusion and future work}\label{sec:7}
In this contribution, we presented a new algorithm to speed up
PDE-constrained optimization problems occurring in the context of enhanced
oil recovery. The algorithm is based on adaptively constructed surrogate
models that make use of deep neural networks for approximating the objective
functional. In each outer iteration of the algorithm, a new surrogate model
is trained with data consisting of full-order function evaluations around the
current control point. Afterwards, an ensemble-based optimization algorithm is
applied to the surrogate to obtain a candidate for the next iteration.
We perform full order model evaluations to validate whether the resulting controls correspond to a local optimum of
the true objective functional. These
function evaluations also serve as training data for constructing the
next surrogate.
\par
Our numerical results confirm that the described algorithm can
accelerate the solution of the enhanced oil recovery optimization problem.
At the same time, in our numerical experiments, the procedure
produces controls with even larger objective function values than those
obtained using only costly full-order model evaluations. However, we
should emphasize that such an improvement in the objective function value is
not guaranteed and, in our case, results from the multi-modal structure of
the objective functional.
\par
The investigated five-spot benchmark problem served as a proof of concept
for our \ROMEnOpt\ algorithm, where FOM evaluations were relatively quickly
accessible, and the input dimension was of moderate size. Future research is
thus devoted to more involved numerical experiments with more significant
complexity.
\par
As indicated in the optimization problem description, we
focused on a scenario with fixed geological properties. However, in practical
applications, these geological parameters are usually unknown and typically
treated by ensemble-based methods, where the ensemble is to be understood not
only with respect to perturbations of the controls for approximating the
gradient but also with respect to different samples of geological properties.
One of the central future research perspectives is incorporating such geological
uncertainty in our algorithm. The main challenge is the high
dimension of the space of possible geological parameters. Naively using these
parameters as additional inputs for the neural network is thus not feasible.
Future research might consider reducing the dimension of the space of
geological parameters by incorporating additional information on the
distribution of such parameters and passing the reduced variables to the
neural networks.
\par
Furthermore, replacing neural networks as surrogate models for the objective
function, for instance, by polynomial approximations obtained via linear
regression or by different machine learning approaches, such as kernel
methods~\cite{hofmann2008}, could be investigated further. The
\ROMEnOpt\ algorithm is formulated in such a way that replacing the surrogate model and
its training is readily possible. Any approximation of the objective
function built from evaluations of the true objective function is feasible
and can directly be used in the algorithm. In addition, the inner iterations
are not restricted to the \EnOpt\ procedure but can also be performed using
different optimization routines. However, we should emphasize that in the
current formulation, no information on the exact gradient, neither of the
true objective functional nor the surrogate model, is required. This might
change when employing different optimization routines. Moreover, the
presented approach is not restricted to the NPV objective functional in
enhanced oil recovery but can be generalized to any scalar-valued quantity
of interest. The algorithm might be of particular relevance in cases where no
direct access to the underlying PDE solutions is possible, and no error
estimation for the surrogate model is available.

\section*{Acknowledgement}
\begin{itemize}
\item Tim Keil, Hendrik Kleikamp, Micheal Oguntola and Mario Ohlberger acknowledge funding by the Deutsche Forschungsgemeinschaft (DFG, German Research Foundation) under Germany's Excellence Strategy EXC 2044 –390685587, Mathematics Münster: Dynamics–Geometry–Structure.
\item Tim Keil and Mario Ohlberger acknowledge funding by the Deutsche Forschungsgemeinschaft under contract OH 98/11-1.
\item Micheal Oguntola and Rolf Lorentzen acknowledge funding from the Research Council of Norway and the industry partners, ConocoPhillips Skandinavia AS, Aker BP ASA, Vår Energi AS, Equinor Energy AS, Neptune Energy Norge AS, Lundin Energy Norway AS, Halliburton AS, Schlumberger Norge AS, and Wintershall Dea Norge AS, of The National IOR Centre of Norway.
\end{itemize}

{\small
	\bibliographystyle{plain}

}

\begin{thebibliography}{10}

\bibitem{abidin2012polymers}
A.~Abidin, T.~Puspasari, and W.~Nugroho.
\newblock Polymers for enhanced oil recovery technology.
\newblock {\em Procedia Chemistry}, 4:11--16, 2012.

\bibitem{ahmadi2015}
M.~A. Ahmadi.
\newblock Developing a robust surrogate model of chemical flooding based on the
  artificial neural network for enhanced oil recovery implications.
\newblock {\em Mathematical Problems in Engineering}, 2015:9, 2015.

\bibitem{banholzer2020adaptive}
S.~Banholzer, T.~Keil, L.~Mechelli, M.~Ohlberger, F.~Schindler, and
  S.~Volkwein.
\newblock An adaptive projected newton non-conforming dual approach for
  trust-region reduced basis approximation of pde-constrained parameter
  optimization, 2020.
\newblock {arXiv}-eprint:2012.11653.

\bibitem{bao2017fully}
K.~Bao, K.-A. Lie, O.~M{\o}yner, and M.~Liu.
\newblock Fully implicit simulation of polymer flooding with mrst.
\newblock {\em Computational Geosciences}, 21(5):1219--1244, 2017.

\bibitem{opm2021}
D.~Baxendale, A.~F. Rasmussen, A.~B. Rustad, T.~Skille, and T.~H. Sandve.
\newblock Opm flow documentation manual manual.
\newblock {\em Open Porous Media Initiative}, 2021.

\bibitem{MR3701994}
P.~Benner, M.~Ohlberger, A.~Patera, G.~Rozza, and K.~Urban, editors.
\newblock {\em Model reduction of parametrized systems}, volume~17 of {\em
  MS\&A. Modeling, Simulation and Applications}.
\newblock Springer, Cham, 2017.
\newblock Selected papers from the 3rd MoRePaS Conference held at the
  International School for Advanced Studies (SISSA), Trieste, October 13--16,
  2015.

\bibitem{bottou2018}
L.~Bottou, F.~E. Curtis, and J.~Nocedal.
\newblock Optimization methods for large-scale machine learning.
\newblock {\em SIAM Rev.}, 60(2):223--311, 2018.

\bibitem{chen2009efficient}
Y.~Chen, D.~S. Oliver, and D.~Zhang.
\newblock Efficient ensemble-based closed-loop production optimization.
\newblock {\em SPE Journal}, 14(04):634--645, 2009.

\bibitem{chen2007reservoir}
Z.~Chen.
\newblock {\em Reservoir simulation: mathematical techniques in oil recovery}.
\newblock SIAM, 2007.

\bibitem{cheraghi2021}
Y.~Cheraghi, S.~Kord, and V.~Mashayekhizadeh.
\newblock Application of machine learning techniques for selecting the most
  suitable enhanced oil recovery method; challenges and opportunities.
\newblock {\em Journal of Petroleum Science and Engineering}, 205:108761, 2021.

\bibitem{elbraechter2018}
D.~Elbrächter, P.~Grohs, A.~Jentzen, and C.~Schwab.
\newblock Dnn expression rate analysis of high-dimensional pdes: Application to
  option pricing.
\newblock {\em Constructive Approximation}, 2021.

\bibitem{fonseca2015quantification}
R.~Fonseca, S.~Kahrobaei, L.~Van~Gastel, O.~Leeuwenburgh, and J.~Jansen.
\newblock Quantification of the impact of ensemble size on the quality of an
  ensemble gradient using principles of hypothesis testing.
\newblock In {\em SPE Reservoir Simulation Symposium}. OnePetro, 2015.

\bibitem{fonseca2017stochastic}
R.~R.-M. Fonseca, B.~Chen, J.~D. Jansen, and A.~Reynolds.
\newblock A stochastic simplex approximate gradient (stosag) for optimization
  under uncertainty.
\newblock {\em International Journal for Numerical Methods in Engineering},
  109(13):1756--1776, 2017.

\bibitem{GHI+2021}
P.~Gavrilenko, B.~Haasdonk, O.~Iliev, M.~Ohlberger, F.~Schindler, P.~Toktaliev,
  T.~Wenzel, and M.~Youssef.
\newblock A full order, reduced order and machine learning model pipeline for
  efficient prediction of reactive flows.
\newblock 2021.
\newblock {arXiv}-eprint:2104.02800.

\bibitem{golzari2015}
A.~Golzari, M.~{Haghighat Sefat}, and S.~Jamshidi.
\newblock Development of an adaptive surrogate model for production
  optimization.
\newblock {\em Journal of Petroleum Science and Engineering}, 133:677--688,
  2015.

\bibitem{gudina2015biosurfactant}
E.~J. Gudi{\~n}a, E.~C. Fernandes, A.~I. Rodrigues, J.~A. Teixeira, and L.~R.
  Rodrigues.
\newblock Biosurfactant production by bacillus subtilis using corn steep liquor
  as culture medium.
\newblock {\em Frontiers in microbiology}, 6:59, 2015.

\bibitem{HOS2021}
B.~Haasdonk, M.~Ohlberger, and F.~Schindler.
\newblock An adaptive model hierarchy for data-augmented training of kernel
  models for reactive flow.
\newblock 2021.
\newblock {arXiv}-eprint:2110.12388.

\bibitem{hastie2009}
T.~Hastie, R.~Tibshirani, and J.~Friedman.
\newblock {\em The Elements of Statistical Learning}.
\newblock 10.1007/978-0-387-84858-7, 2009.

\bibitem{he2015}
K.~{He}, X.~{Zhang}, S.~{Ren}, and J.~{Sun}.
\newblock Delving deep into rectifiers: Surpassing human-level performance on
  imagenet classification.
\newblock In {\em 2015 IEEE International Conference on Computer Vision
  (ICCV)}, pages 1026--1034, 2015.

\bibitem{hofmann2008}
T.~Hofmann, B.~Schölkopf, and A.~J. Smola.
\newblock {Kernel methods in machine learning}.
\newblock {\em The Annals of Statistics}, 36(3):1171 -- 1220, 2008.

\bibitem{holmes1983enhancements}
J.~Holmes.
\newblock Enhancements to the strongly coupled, fully implicit well model:
  wellbore crossflow modeling and collective well control.
\newblock In {\em SPE Reservoir Simulation Symposium}. OnePetro, 1983.

\bibitem{holmes1998application}
J.~Holmes, T.~Barkve, and O.~Lund.
\newblock Application of a multisegment well model to simulate flow in advanced
  wells.
\newblock In {\em European petroleum conference}. OnePetro, 1998.

\bibitem{MR4269464}
T.~Keil, L.~Mechelli, M.~Ohlberger, F.~Schindler, and S.~Volkwein.
\newblock A non-conforming dual approach for adaptive trust-region reduced
  basis approximation of {PDE}-constrained parameter optimization.
\newblock {\em ESAIM Math. Model. Numer. Anal.}, 55(3):1239--1269, 2021.

\bibitem{lecun2015}
Y.~LeCun, Y.~Bengio, and G.~Hinton.
\newblock Deep learning.
\newblock {\em Nature}, 521:436--44, 2015.

\bibitem{lee2011}
J.-Y. Lee, H.-J. Shin, and J.-S. Lim.
\newblock Selection and evaluation of enhanced oil recovery method using
  artificial neural network.
\newblock {\em Geosystem Engineering}, 14:157 -- 164, 2011.

\bibitem{liu1989}
D.~C. Liu and J.~Nocedal.
\newblock On the limited memory {BFGS} method for large scale optimization.
\newblock {\em Mathematical Programming}, 45:503--528, 1989.

\bibitem{lu2020joint}
R.~Lu and A.~Reynolds.
\newblock Joint optimization of well locations, types, drilling order, and
  controls given a set of potential drilling paths.
\newblock {\em SPE Journal}, 25(03):1285--1306, 2020.

\bibitem{lye2021}
K.~O. Lye, S.~Mishra, D.~Ray, and P.~Chandrashekar.
\newblock Iterative surrogate model optimization (ismo): An active learning
  algorithm for pde constrained optimization with deep neural networks.
\newblock {\em Computer Methods in Applied Mechanics and Engineering},
  374:113575, 2021.

\bibitem{milk2016}
R.~Milk, S.~Rave, and F.~Schindler.
\newblock {pyMOR} -- generic algorithms and interfaces for model order
  reduction.
\newblock {\em SIAM J. Sci. Comput.}, 38(5):S194--S216, jan 2016.

\bibitem{montgomery2015introduction}
D.~C. Montgomery, C.~L. Jennings, and M.~Kulahci.
\newblock {\em Introduction to time series analysis and forecasting}.
\newblock John Wiley \& Sons, 2015.

\bibitem{nocedal2006numerical}
J.~Nocedal and S.~Wright.
\newblock {\em Numerical optimization}.
\newblock Springer Science \& Business Media, 2006.

\bibitem{oguntola2020robust}
M.~B. Oguntola and R.~J. Lorentzen.
\newblock On the robust value quantification of polymer eor injection
  strategies for better decision making.
\newblock In {\em ECMOR XVII}, volume 2020, pages 1--25. European Association
  of Geoscientists \& Engineers, 2020.

\bibitem{oguntola2021ensemble}
M.~B. Oguntola and R.~J. Lorentzen.
\newblock Ensemble-based constrained optimization using an exterior penalty
  method.
\newblock {\em Journal of Petroleum Science and Engineering}, 207:109165, 2021.

\bibitem{pancholi2020experimental}
S.~Pancholi, G.~S. Negi, J.~R. Agarwal, A.~Bera, and M.~Shah.
\newblock Experimental and simulation studies for optimization of
  water--alternating-gas (co2) flooding for enhanced oil recovery.
\newblock {\em Petroleum Research}, 5(3):227--234, 2020.

\bibitem{paszke2019}
A.~Paszke, S.~Gross, F.~Massa, A.~Lerer, J.~Bradbury, G.~Chanan, T.~Killeen,
  Z.~Lin, N.~Gimelshein, L.~Antiga, A.~Desmaison, A.~Kopf, E.~Yang, Z.~DeVito,
  M.~Raison, A.~Tejani, S.~Chilamkurthy, B.~Steiner, L.~Fang, J.~Bai, and
  S.~Chintala.
\newblock Pytorch: An imperative style, high-performance deep learning library.
\newblock In H.~Wallach, H.~Larochelle, A.~Beygelzimer, F.~d'~Alch\'{e}-Buc,
  E.~Fox, and R.~Garnett, editors, {\em Advances in Neural Information
  Processing Systems 32}, pages 8024--8035. Curran Associates, Inc., 2019.

\bibitem{petersen2018}
P.~Petersen and F.~Voigtlaender.
\newblock Optimal approximation of piecewise smooth functions using deep relu
  neural networks.
\newblock {\em Neural Networks}, 108:296 -- 330, 2018.

\bibitem{prechelt1997}
L.~Prechelt.
\newblock Early stopping - but when?
\newblock In {\em Neural Networks: Tricks of the Trade, volume 1524 of LNCS,
  chapter 2}, pages 55--69. Springer-Verlag, 1997.

\bibitem{flo2019open}
A.~F. Rasmussen, T.~H. Sandve, K.~Bao, A.~Lauser, J.~Hove, B.~Skaflestad,
  R.~Kl{\"o}fkorn, M.~Blatt, A.~B. Rustad, O.~S{\ae}vareid, et~al.
\newblock The open porous media flow reservoir simulator.
\newblock {\em Computers \& Mathematics with Applications}, 81:159--185, 2021.

\bibitem{rumelhart1986}
D.~E. Rumelhart, G.~E. Hintont, and R.~J. Williams.
\newblock {Learning representations by back-propagating errors}.
\newblock {\em Nature}, 323(6088):533--536, 1986.

\bibitem{hossein2021}
H.~Saberi, E.~Esmaeilnezhad, and H.~J. Choi.
\newblock Artificial neural network to forecast enhanced oil recovery using
  hydrolyzed polyacrylamide in sandstone and carbonate reservoirs.
\newblock {\em Polymers}, 13(16), 2021.

\bibitem{sarma2006efficient}
P.~Sarma, L.~J. Durlofsky, K.~Aziz, and W.~H. Chen.
\newblock Efficient real-time reservoir management using adjoint-based optimal
  control and model updating.
\newblock {\em Computational Geosciences}, 10(1):3--36, 2006.

\bibitem{stordal2016theoretical}
A.~S. Stordal, S.~P. Szklarz, and O.~Leeuwenburgh.
\newblock A theoretical look at ensemble-based optimization in reservoir
  management.
\newblock {\em Mathematical Geosciences}, 48(4):399--417, 2016.

\bibitem{van2016well}
S.~L. Van and B.~H. Chon.
\newblock Well-pattern investigation and selection by surfactant-polymer
  flooding performance in heterogeneous reservoir consisting of interbedded
  low-permeability layer.
\newblock {\em Korean Journal of Chemical Engineering}, 33(12):3456--3464,
  2016.

\bibitem{wang2008key}
D.~Wang, R.~S. Seright, Z.~Shao, J.~Wang, et~al.
\newblock Key aspects of project design for polymer flooding at the daqing
  oilfield.
\newblock {\em SPE Reservoir Evaluation \& Engineering}, 11(06):1--117, 2008.

\bibitem{wang2021fast}
S.~Wang, M.~A. Bhouri, and P.~Perdikaris.
\newblock Fast pde-constrained optimization via self-supervised operator
  learning, 2021.
\newblock {arXiv}-eprint:2110.13297.

\bibitem{wolfe1969}
P.~Wolfe.
\newblock Convergence conditions for ascent methods.
\newblock {\em SIAM Review}, 11(2):226--235, 1969.

\bibitem{wolfe1971}
P.~Wolfe.
\newblock Convergence conditions for ascent methods. ii: Some corrections.
\newblock {\em SIAM Review}, 13(2):185--188, 1971.

\bibitem{xu2018production}
L.~Xu, H.~Zhao, Y.~Li, L.~Cao, X.~Xie, X.~Zhang, and Y.~Li.
\newblock Production optimization of polymer flooding using improved monte
  carlo gradient approximation algorithm with constraints.
\newblock {\em Journal of Circuits, Systems and Computers}, 27(11):1850167,
  2018.

\bibitem{yarotsky2017}
D.~Yarotsky.
\newblock Error bounds for approximations with deep {ReLU} networks.
\newblock {\em Neural Networks}, 94:103 -- 114, 2017.

\bibitem{Zahr2015}
M.~J. Zahr and C.~Farhat.
\newblock Progressive construction of a parametric reduced-order model for
  {PDE}-constrained optimization.
\newblock {\em Int. J. Numer. Meth. Engng}, 102:1111--1135, 2015.

\bibitem{zhang2017well}
Y.~Zhang, R.~Lu, F.~Forouzanfar, and A.~C. Reynolds.
\newblock Well placement and control optimization for wag/sag processes using
  ensemble-based method.
\newblock {\em Computers \& Chemical Engineering}, 101:193--209, 2017.

\bibitem{zhou2013optimal}
K.~Zhou, J.~Hou, X.~Zhang, Q.~Du, X.~Kang, and S.~Jiang.
\newblock Optimal control of polymer flooding based on simultaneous
  perturbation stochastic approximation method guided by finite difference
  gradient.
\newblock {\em Computers \& chemical engineering}, 55:40--49, 2013.

\end{thebibliography}

\end{document}